\documentclass[12pt,a4paper]{article}
\usepackage{amsmath}
\usepackage{graphicx}
\allowdisplaybreaks[1]

\newcommand{\pr}{\rightarrow}

\newcommand{\ba}{\begin{array}}
\newcommand{\ea}{\end{array}}

\newcommand{\vart}{\vartheta}
\newcommand{\varp}{\varphi}
\newcommand{\eps}{\varepsilon}

\newcommand{\il}{\int\limits}
\newenvironment{inspring}[1]%
{\begin{list}{}{\setlength{\rightmargin}{0cm}
                \setlength{\listparindent}{0cm}
                \settowidth{\labelwidth}{\mbox{#1}}
                \setlength{\leftmargin}{1.1\labelwidth}
                \setlength{\labelsep}{.1\labelwidth}}}%
{\end{list}}
\newcommand{\ITEM}[1]{\item[#1\hfill]}
\newcommand{\bi}[1]{\begin{inspring}{#1}}
\newcommand{\ei}{\end{inspring}}

\newcommand{\bfx}{{\bf x}}

\newcommand{\bfn}{{\bf n}}
\newcommand{\bom}{\boldsymbol{\omega}}
\newcommand{\bxi}{\boldsymbol{\xi}}
\newcommand{\bfeta}{\boldsymbol{\eta}}
\newcommand{\beq}{\begin{equation}}
\newcommand{\eq}{\end{equation}}
\catcode`\@=11

\font\tenmsa=msam10 \font\sevenmsa=msam7 \font\fivemsa=msam5
\font\tenmsb=msbm10 \font\sevenmsb=msbm7 \font\fivemsb=msbm5
\newfam\msafam
\newfam\msbfam
\textfont\msafam=\tenmsa  \scriptfont\msafam=\sevenmsa
  \scriptscriptfont\msafam=\fivemsa
\textfont\msbfam=\tenmsb  \scriptfont\msbfam=\sevenmsb
  \scriptscriptfont\msbfam=\fivemsb

\def\Bbb{\ifmmode\let\next\Bbb@\else
 \def\next{\errmessage{Use \string\Bbb\space only in math mode}}\fi\next}
\def\Bbb@#1{{\Bbb@@{#1}}}
\def\Bbb@@#1{\fam\msbfam#1}
\newcommand{\dR}{{\Bbb R}}

\numberwithin{thm}{section}

\title{Generalized 3D Zernike functions for analytic construction of band-limited line-detecting wavelets}
\author{A.J.E.M.\ Janssen \\
Eindhoven University of Technology, \\
Department of Mathematics and Computer Science, \\
P.O.\ Box 513, 5600 MB Eindhoven, The Netherlands. \\
\{a.j.e.m.janssen@tue.nl\}}
\date{}

\begin{document}
\maketitle
\mbox{} \\ \\ \\ \\ \\
\noindent
{\bf Abstract.} \\
We consider 3D versions of the Zernike polynomials that are commonly used in 2D in optics and lithography. We generalize the 3D Zernike polynomials to functions that vanish to a prescribed degree $\alpha\geq0$ at the rim of their supporting ball $\rho\leq1$. The analytic theory of the 3D generalized Zernike functions is developed, with attention for computational results for their Fourier transform, Funk and Radon transform, and scaling operations. The Fourier transform of generalized 3D Zernike functions shows less oscillatory behaviour and more rapid decay at infinity, compared to the standard case $\alpha=0$, when the smoothness parameter $\alpha$ is increased beyond 0. The 3D generalized Zernike functions can be used to expand smooth functions, supported by the unit ball and vanishing at the rim and the origin of the unit ball, whose radial and angular dependence is separated. Particular instances of the latter functions (prewavelets) yield, via the Funk transform and the Fourier transform, an anisotropic function that can be used for a band-limited line-detecting wavelet transform, appropriate for analysis of 3D medical data containing elongated structures. We present instances of prewavelets, with relevant radial functions, that allow analytic computation of Funk and Fourier transform. A key step here is to identify the special form that is assumed by the expansion coefficients of a separable function on the unit ball with respect to generalized 3D Zernike functions. A further issue is how to scale a function on the unit ball while maintaining its supporting set, and this issue is solved in a particular form.

\section{Introduction} \label{sec1}
\subsection{Motivation} \label{subsec1.1}
\mbox{} \\[-9mm]

In the signal processing of 2D or 3D medical band-limited image data containing elongated structures, a special form of the wavelet transform is used. This special wavelet transform is built from an anisotropic wavelet $\psi$ whose modulus has a contour plot with relatively long ridges in one particular direction. The Fourier transform $F={\cal F}\psi$ of $\psi$ is concentrated near the line (2D) or the plane (3D) perpendicular to this preferred direction, as long as contained in the supporting disk (2D) or ball (3D) of $F={\cal F}\psi$ that we assume to be the unit disk or ball. In the case of 2D image processing, one may choose a function $F$ that factorizes as $AB$, where $A$ depends on the angular variable $\varp$ and $B$ on the radial variable $\rho$. Here $A$ is even around $\pi/2$ and non-vanishing in a relatively small interval around $\pi/2$, and $B$ is positive and away from 0 on a relatively large or small interval contained in $(0,1)$, depending on whether one wants an all-scale or a multi-scale wavelet transform, where in the latter case variable-scaled versions of $B$ occur. The inverse Fourier transform $\psi={\cal F}^{-1}F$ is then disk-limited (i.e., a function with Fourier transform supported by the unit disk), and the contour plot of $|\psi|$ has ridges mainly parallel to the $x$-axis. In the case of 3D image processing, one still departs from an $F$ of the factorized form $AB$, with $A$ depending on the spherical angular variables $\vart=[0,\pi]$ and $\varp\in[0,2\pi]$. This $A(\vart,\varp)$ is often of the form $h(\vart)$ in which $h$ is a smooth function supported by a relatively small interval $[0,\vart_0]\subset[0,\pi]$. Next, $A$ is replaced by its Funk transform. The Funk transform maps the function $A$, considered as a function defined on the unit sphere, onto the function $A_{{\rm Funk}}$ defined on the unit sphere whose value $A_{{\rm Funk}}(\bfn)$ at any unit vector $\bfn$ is given by the average of $A$ over the great circle perpendicular to $\bfn$. For the $A$ at hand, the Funk transform $A_{{\rm Funk}}$ is rotationally symmetric about the vertical axis $(\vart=0)$ and has as a supporting set the meridional zone $|\vart-\pi/2|\leq\vart_0$ near the plane perpendicular to the vertical axis.

In the proposals considered until now \cite{ref1}--\cite{ref5}, the choice of the radial function $B$ has been defined in terms of $B$-splines, with logarithmic dependence on $\rho$, or in terms of Gaussians. In the first proposal, the computation of $\psi$ as a Fourier transform is analytically not feasible, and one has to resort to DFT-methods with the usual sampling issues. In the second proposal, one has a $B$ that has, strictly spoken, an unbounded support. In addition, in both cases, special measures have to be taken in order that $\psi$ and its scaled versions can be considered to be disk or ball-limited while having sufficient decay and a controlled amount of oscillations.

In this report, we propose to use orthogonal functions $Z$ on the unit disk or ball to expand separable functions $F$, with the requirements
\bi{a.0}
\ITEM{a.} $Z$ decays at the rim of the unit disk or ball,
\ITEM{b.} the Fourier transform of $Z$ has an analytic form,
\ITEM{c.} angular and radial dependence of $Z$ is separated,
\ITEM{d.} the expansion in $Z$'s is feasible for a set of separable functions $F$ containing cases relevant in the present context,
\ITEM{e.} the scaling issue is resolved.
\ei
For the 2D case, a family of orthogonal systems satisfying the requirements a, b, c has been considered in \cite{ref6}, viz.\ the generalized Zernike functions $Z_n^{m,\alpha}$. For $\alpha>-1$ and integer $n$ and $m$ such that $n-|m|$ is even and non-negative $Z_n^{m,\alpha}$ is given by
\beq \label{e1}
Z_n^m(\nu,\mu)\equiv Z_n^m(\rho,\varp)=R_n^{|m|,\alpha}(\rho)\,e^{im\varp}~,~~~~~~0\leq\rho\leq1\,,~~0\leq\rho\leq 2\pi~,
\eq
where we write $\nu+i\mu$ with $\nu,\mu\in\dR$ and $\nu^2+\mu^2\leq1$ as $\rho\,\exp(i\varp)$ with $0\leq\rho\leq1$, $0\leq\varp\leq2\pi$, and where
\beq \label{e2}
R_n^{|m|}(\rho)=\left\{\ba{llll}
(1-\rho^2)^{\alpha}\,\rho^{|m|}\,P_{\frac{n-|m|}{2}}^{(\alpha,|m|)}(2\rho^2-1) & \!\!, & ~~~0\leq\rho\leq1 & \!\!, \\[3mm]
0 & \!\!, & ~~~\rho>1 & \!\!,
\ea\right.
\eq
and $P_k^{(\alpha,\beta)}(x)$ is the Jacobi polynomial corresponding to the weight function $(1-x)^{\alpha}(1+x)^{\beta}$, $-1\leq x\leq1$, of degree $k$. The functions $Z_n^{m,\alpha}$ have been considered earlier, in a general setting of unit balls in $\dR^d$, in the context of the Radon transform in \cite{ref7}. The case that $\alpha=0$ yields the standard Zernike circle polynomials that find wide-spread application in optics, lithography, acoustics \cite{ref8}--\cite{ref10}.

For the 3D case, we consider in this report the generalized 3D Zernike functions. These are defined for $\alpha>{-}1$ and non-negative integer $n,l$ such that $n-l$ is even and non-negative and $m={-}l,{-}l+1,...,l$ by
\beq \label{e3}
Z_{nl}^{m,\alpha}(\nu,\mu,\sigma)=R_n^{l,\alpha}(\rho)\,Y_l^m(\vart,\varp)~,~~~~~~0\leq\rho\leq1\,,~~0\leq\vart\leq\pi\,,~~0\leq\varp\leq2\pi~,
\eq
where we employ spherical coordinates
\beq \label{e4}
(\nu,\mu,\sigma)=\bom=(\rho\sin\vart\cos\varp,\rho\sin\vart\sin\varp,\rho\cos\vart)~.
\eq
The radial function $R_n^{l,\alpha}(\rho)$ in (\ref{e3}) is given by
\beq \label{e5}
R_n^{l,\alpha}(\rho)=\rho^l(1-\rho^2)^{\alpha}\,P_{\frac{n-l}{2}}^{(\alpha,l+1/2)}(2\rho^2-1)~,~~~~~~0\leq\rho\leq1~,
\eq
and vanishes for $\rho>1$, and the angular function $Y_l^m(\vart,\varp)$ is the spherical harmonic
\beq \label{e6}
Y_l^m(\vart,\varp)=({-}1)^m\Bigl(\frac{2l+1}{4\pi}\:\frac{(l-|m|!}{(l+m)!}\Bigr)^{1/2}\,P_l^{|m|}(\cos\vart)\,e^{im\varp}~,
\eq
normalized to have unit $L^2$-norm on the 3D sphere. Note that we use unnormalized angular functions $\exp(im\varp)$ in the definition of the 2D generalized Zernike functions in (\ref{e1}) for reasons of consistency with \cite{ref6}. The case $\alpha=0$ in (\ref{e3}) yield the standard Zernike ball polynomials that are considered, with yet another normalization convention, in \cite{ref11}--\cite{ref12}.

A great deal of the basic and more advanced properties of the generalized 2D Zernike functions have been developed in \cite{ref6}, but the above issues d and e have not been addressed there. In the present report, the emphasis is on developing the theory and results for the 3D case in which all 5 issues a--e are addressed. In this effort, we shall present the 2D versions of the results, when not already covered by \cite{ref6}, by appending an additional superscript 2 to $Z$, $R$, etc., to distinguish from the 3D case.

Since the radial and angular dependence are separated in the generalized 2D and 3D Zernike functions, the expansion coefficients of a factorized $F=AB$ factorize to a large extent as well. However, the radial functions in (\ref{e2}) and (\ref{e5}) contain also the angular index, $m$ and $l$, respectively, yielding for any $m$ and $l$ a different expansion of $B$ into radial functions. The set of relevant angular orders $m$ or $l$ that occur, is determined by smoothness of the angular factor $A$. In order to ensure that the coefficients in the $B$-expansions remain tractable in size when $m$ or $l$ varies, the degree to which $B$ vanishes at $\rho=0$ should be chosen sufficiently high.

In 3D, the Funk transform retains only spherical harmonics $Y_l^m$ with even value of $l$, and so we only need the expansion of $B$ into radial functions $R_n^{l,\alpha}$ with $l=0,2,...\,$. Similarly, in the 2D case, the function $A$ is even around $\varp=\pi/2$ and this implies that only the expansions of $B$ into radial functions $R_n^{|m|,\alpha}$ with even $m$ are required.

\subsection{Overview of the results} \label{subsec1.2}
\mbox{} \\[-9mm]

In Sec.~\ref{sec2} we present the more basic properties of the generalized 3D Zernike functions. Thus we consider orthogonality and normalization of the $Z$'s, we show how separability of an $F=AB$ is reflected by separability of its expansion coefficients, and we compute the Fourier transforms of the $Z$'s using the Funk-Hecke formula. In Sec.~\ref{sec3}, we present the required results concerning scaling 3D generalized Zernike functions. In Sec.~\ref{sec4}, we present an explicit result for the Radon transform of 3D generalized Zernike functions and an integral representation, in terms of Gegenbauer polynomials, of the corresponding radial functions that follows from this Radon transform result. This integral representation can be used to prove a recursion of the Shakibaei-Paramesran type \cite{ref13} to compute all radial functions at a particular $\rho\in[0,1]$. The next few sections are devoted to the study and computation of the expansion coefficients of various angular order $l$ of a radial function $B$ occurring as radial factor of a separable $F$. It is an important fact, to be proved in Sec.~\ref{sec5}, that the coefficients required for $l+2$ can be expressed explicitly and in finite terms in those required for $l$. Since we may restrict to even $l$, it thus follows that all required expansions of $B$ can be obtained from the one with angular order $l=0$. In Sec.~\ref{sec6} we consider radial functions $B(\rho)$ of the form
\beq \label{e7}
B(\rho)=\rho^{\beta}\,C(\eps\rho)(1-\rho^2)^{\delta}~,
\eq
where $\beta$ and $\delta$ are non-negative integers with $\beta$ even, $\eps\in[0,1]$ is a scaling parameter, and $C(\rho)$ is a smooth function defined for $\rho\in[0,1]$ and not necessarily vanishing at $\rho=0$ or 1. We describe a procedure how to get from the expansion of $C(\rho)$ into radial functions $R_{2k}^{0,0}$ all such expansions for $B(\rho)$ in (\ref{e7}). The choice
\beq \label{e8}
C(\rho)=(1-\rho^2)^{\eta}
\eq
in (\ref{e7}) with integer $\eta\geq0$ is considered in Sec.~\ref{sec7} as a special case of radial functions that admit explicit expansion into all or particular radial functions $R_n^{l,\alpha}(\rho)$. Finally, in Sec.~\ref{sec8} we present an example of a prewavelet radial function meant for an all-scale wavelet transform, and in Sec.~\ref{sec9} we present such an example, of the type (\ref{e7}), for a multi-scale wavelet transform. For these two cases, there are explicit expansions for all required angular orders $l$. Hence, when we consider the expansion coefficients of the angular functions $A$ in spherical harmonics $Y_l^m$ as being given, this yields an explicit, analytic result for the wavelet $\psi={\cal F}F$.

\section{Basic properties of 3D generalized Zernike functions} \label{sec2}
\mbox{} \\[-9mm]

We present basic properties of the 3D generalized Zernike functions, defined in (\ref{e3}--\ref{e6}) in terms of Jacobi polynomials and spherical harmonics. Relevant reference for the latter polynomials and special functions are \cite{ref14}, Ch.~18 and \S 14.30, \cite{ref15}, Ch.~4, and \cite{ref16}, Ch.~5.

\subsection{Orthogonality} \label{subsec2.1}
\mbox{} \\[-9mm]

We have for integer $n_1,n_2,l_1,l_2,m_1,m_2$ such that $n_1,n_2,l_1,l_2\geq0$, $n_1-l_1$ and $n_2-l_2$ even and non-negative, $|m_1|\leq l_1$ and $|m_2|\leq l_2$, and $\alpha>{-}1$
\begin{eqnarray} \label{e9}
& \mbox{} & \iiint\limits_{\nu^2+\mu^2+\sigma^2\leq1}
\,Z_{n_1l_1}^{m_1,\alpha}(\nu,\mu,\sigma)(Z_{n_2l_2}^{m_2,\alpha}(\nu,\mu,\sigma))^{\ast} \,\frac{d\nu d\mu d\sigma}{(1-\rho^2)^{\alpha}} \nonumber \\[3.5mm]
& & =~\il_0^1\,R_{n_1}^{l_1,\alpha}(\rho)\,R_{n_2}^{l_2,\alpha}(\rho)\,\frac{\rho^2\,d\rho}{(1-\rho^2)^{\alpha}}\cdot \il_0^{\pi}\il_0^{2\pi}\,Y_{l_1}^{m_1}(\vart,\varp)(Y_{l_2}^{m_2}(\vart,\varp))^{\ast} \nonumber \\[2mm]
& & \hspace*{8cm}\cdot\:\sin\vart\,d\vart d\varp~= \nonumber \\[2mm]
& & =~N_{nl}^{\alpha}\,\delta_{n_1n_2}\cdot\delta_{m_1m_2}\,\delta_{l_1l_2}~,
\end{eqnarray}
where $\delta$ is Kronecker's delta. Furthermore, for $n_1=n_2=n=l+2p$ with $l,p=0,1,...\,$, we have
\beq \label{e10}
N_{nl}^{\alpha}=\frac{1}{2(n+\alpha+3/2)}~\frac{(p+1)_{\alpha}}{(p+l+3/2)_{\alpha}}~,
\eq
where we employ the generalized Pochhammer symbol
\beq \label{e11}
(x)_{\alpha}=\frac{\Gamma(x+\alpha)}{\Gamma(x)}~.
\eq
Similarly, in the 2D case, we have for integer $n_1,m_1,n_2,m_2$ with $n_1-|m_1|$ and $n_2-|m_2|$ even and non-negative, and $\alpha>{-}1$
\beq \label{e12}
\iint\limits_{\nu^2+\mu^2\leq1}{}^2Z_{n_1}^{m_1,\alpha}(\nu,\mu)(^2Z_{n_2}^{m_2,\alpha}(\nu,\mu))^{\ast}\, \frac{d\nu d\mu}{(1-\rho^2)^{\alpha}}=2\pi\,{}^2N_{nm}^{\alpha}\,\delta_{n_1n_2}\,\delta_{m_1m_2}~,
\eq
where
\beq \label{e13}
{}^2N_{nm}^{\alpha}=\il_0^1\,|R_n^{|m|,\alpha}(\rho)|^2\,\frac{\rho\,d\rho}{(1-\rho^2)^{\alpha}}=\frac{1}{2(n+\alpha+1)}\, \frac{(p+1)_{\alpha}}{(p+|m|+1)_{\alpha}}~,
\eq
with $n=|m|+2p$ and $p=0,1,...\,$.

\subsection{Expansion coefficients of separable functions} \label{subsec2.2}
\mbox{} \\[-9mm]

For
\beq \label{e14}
F(\nu,\mu,\sigma)\equiv F(\vart,\varp,\rho)=A(\vart,\varp)\,B(\rho)~,
\eq
with $A(\vart,\varp)$ and $B(\rho)$ smooth functions of $(\vart,\varp)\in[0,\pi]\times[0,2\pi]$ and $\rho\in[0,1]$, respectively, we have
\beq \label{e15}
F=\sum_{n,m,l}\,c_{nl}^{m,\alpha}\,Z_{nl}^{m,\alpha}~.
\eq
Here, summation is over all integer $n,m,l$ with $n,l\geq0$ and $n-l$ even and non-negative and $|m|\leq l$, and
\begin{eqnarray} \label{e16}
c_{nl}^{m,\alpha} & = & \frac{1}{N_{nl}^{\alpha}}~~\iiint\limits_{\nu^2+\mu^2+\sigma^2\leq1}\, F(\nu,\mu,\sigma)(Z_{nl}^{m,\alpha}(\nu,\mu,\sigma))^{\ast}\,\frac{d\nu d\mu d\sigma}{(1-\rho^2)^{\alpha}} \nonumber \\[3.5mm]
& = & a_l^m(A)\,\frac{1}{N_{nl}^{\alpha}}\,b_n^{l,\alpha}(B)~,
\end{eqnarray}
with $a_l^m(A)$ and $b_n^{l,\alpha}(B)$ given by
\beq \label{e17}
a_l^m(A)=\il_0^{\pi}\il_0^{2\pi}\,A(\vart,\varp)\,(Y_l^m(\vart,\varp))^{\ast}\sin\vart\,d\vart d\varp~,
\eq
\beq \label{e18}
b_n^{l,\alpha}(B)=\il_0^1\,B(\rho)\,R_n^{l,\alpha}(\rho)\,\frac{\rho^2\,d\rho}{(1-\rho^2)^{\alpha}}~.
\eq
Thus the $c$'s factorize as $a_l^m(A)\,(N_{nl}^{\alpha})^{-1}\,b_n^{l,\alpha}(B)$, where $a$ depends only on $A$ and $b$ depends only on $B$. The angular index $l$ is present in both factors $a$ and $b$. Thus, there are the expansions
\beq \label{e19}
A(\vart,\varp)=\sum_{l=0}^{\infty}\,\sum_{m={-}l}^l\,a_l^m(A)\,Y_l^m(\vart,\varp)~,
\eq
and, for all $l=0,1,...\,$,
\beq \label{e20}
B(\rho)=\sum_{n=l,l+2,...}\,
\frac{b_n^{l,\alpha}}{N_{nl}^{\alpha}}\,R_n^{l,\alpha}(\rho)~.
\eq

Similarly, we have in the 2D case for $F(\nu,\mu)\equiv F(\varp,\rho)=A(\varp)\,B(\rho)$ that
\beq \label{e21}
F=\sum_{n,m}\,c_n^{m,\alpha}\,{}^2Z_n^m~,
\eq
with summation over all integer $n,m$ such that $n-|m|$ is even and non-negative. The $c$'s are given as
\beq \label{e22}
c_n^{m,\alpha}=a^m(A)\,\frac{1}{{}^2N_{nm}^{\alpha}}\,b_n^{m,\alpha}(B)~,
\eq
with
\beq \label{e23}
a^m(A)=\frac{1}{2\pi}\,\il_0^{2\pi}\,A(\varp)\,e^{-im\varp}\,d\varp~,~~~~~ b_n^{m,\alpha}(B)=\il_0^1\,B(\rho)\,{}^2R_n^{|m|,\alpha}(\rho)\,\frac{\rho\,d\rho}{(1-\rho^2)^{\alpha}}~.
\eq
There are the expansions
\beq \label{e24}
A(\varp)=\sum_{m={-}\infty}^{\infty}\,a^m(A)\,e^{im\varp}~,
\eq
and, for every integer $m$,
\beq \label{e25}
B(\rho)=\sum_{n=|m|,|m|+2,...}\,\frac{b^{m,\alpha}(B)}{{}^2N_{nm}^{\alpha}}\,{}^2R_n^{|m|,\alpha}(\rho)~.
\eq

\subsection{Funk-Hecke formula for spherical harmonics} \label{subsec2.3}
\mbox{} \\[-9mm]

Denote the unit ball in $\dR^3$ by $B$ and the unit sphere in $\dR^3$ by $S$. Let $f$ be integrable over $[{-}1,1]$. Then, see \cite{ref17}, Theorem~1, for $\bom'\in S$ and integer $l,m$ with $|m|\leq l$,
\beq \label{e26}
\il_{\bom\in S}\,f(\bom\cdot\bom')\,Y_l^m(\bom)\,dS=2\pi\lambda_l\,Y_l^m(\bom')~,
\eq
where
\beq \label{e27}
\lambda_l=\il_{-1}^1\,P_l(t)\,f(t)\,dt~,
\eq
with $P_l$ the Legendre polynomial of degree $l$. Specifically, taking $f(t)=\delta_0(t)$ (the Dirac delta function at 0; to be approximated by smooth functions, etc.), we get
\beq \label{e28}
\il_{\bom\in S}\,Y_l^m(\bom)\,\delta_0(\bom\cdot\bom')\,dS=2\pi\, P_l(0)\,Y_l^m(\bom')~,~~~~~~\bom'\in S~.
\eq
The left-hand side of (\ref{e28}) is the integral of $Y_l^m$ over the great circle of radius 1 and center 0 perpendicular to $\bom'$. As to the right-hand side, we have $P_{2k+1}(0)=0$, and
\beq \label{e29}
P_{2k}(0)=({-}1)^k\,\frac{(1/2)_k}{k!}
\eq
for $k=0,1,...\,$.

\subsection{Fourier transform of 3D generalized Zernike functions} \label{subsec2.4}
\mbox{} \\[-9mm]

We have for $\bfx\in\dR^3$
\begin{eqnarray} \label{e30}
{\cal F}\,[Z_{nl}^{m,\alpha}](\bfx)& = & \iiint\limits_{\bom\in B}\,e^{2\pi i\bom\cdot\bfx} \,Z_{nl}^{m,\alpha}(\bom)\,d\bom \nonumber \\[3.5mm]
& = & \il_0^1\il_{\bfeta\in S}\,e^{2\pi i\rho\bfeta\cdot\bfx}\,R_n^{l,\alpha}(\rho)
,Y_l^m(\bfeta)\,d\bfeta\,\rho^2\,d\rho~.
\end{eqnarray}
Write $\bfx=r\bxi$ with $r\geq0$ and $\bxi\in S$, so that
\beq \label{e31}
2\pi i\rho\bfeta\cdot\bfx=is\bfeta\cdot\bxi~,~~~~~~s=2\pi\rho r~.
\eq
We have by the Funk-Hecke formula
\begin{eqnarray} \label{e32}
\iint\limits_{\bfeta\in S}\,e^{is\bfeta\cdot\bxi}\,Y_l^m(\bfeta)\,d\bfeta & = & 2\pi\,Y_l^m(\bxi)\cdot\il_{-1}^1\,e^{ist}\,P_l(t)\,dt \nonumber \\[3.5mm]
& = & 4\pi i^l\,j_l(s)\,Y_l^m(\bxi)~,
\end{eqnarray}
where we have used \cite{ref14}, 18.17.19 on p.~456, with $j_l$ the spherical Bessel function of order $l$. It follows that
\beq \label{e33}
{\cal F}\,[Z_{nl}^{m,\alpha}](\bfx)=4\pi i^l\,Y_l^m(\bxi)\,\il_0^1\,R_n^{l,\alpha}(\rho)\,j_l(2\pi r\rho)\,\rho^2\,d\rho~.
\eq
For the remaining integral, we show below that for $q>0$
\beq \label{e34}
\il_0^1\,R_n^{l,\alpha}(\rho)\,j_l(q\rho)\,\rho^2\,d\rho=({-}1)^p\,2^{\alpha}(p+1)_{\alpha}\, \frac{j_{n+\alpha+1}(q)}{q^{\alpha+1}}~,
\eq
where the right-hand side of (\ref{e34}) equals $\frac12 B(\alpha+1,3/2)\,\delta_{n0}$ at $q=0$, with
\beq \label{e35}
B(a,b)=\frac{\Gamma(a)\,\Gamma(b)}{\Gamma(a+b)}~.
\eq
Using (\ref{e34}) in (\ref{e33}), using $n=l+2p$, we get
\beq \label{e36}
{\cal F}\,[Z_{nl}^{m,\alpha}](\bfx)=2\pi i^n(p+1)_{\alpha}\,\frac{j_{n+\alpha+1}(2\pi|\bfx|)}{(\pi|\bfx|)^{\alpha+1}}\,Y_l^m \Bigl(\frac{\bfx}{|\bfx|}\Bigr)~.
\eq

We now show (\ref{e34}). We use the power series expansion
\beq \label{e37}
j_a(z)=\sqrt{\dfrac{\pi}{2z}}\,J_{a+1/2}(z)=\tfrac12\,\sqrt{\pi}\,(\tfrac12 z)^a\,\sum_{k=0}^{\infty}\,\frac{({-}\tfrac14 z^2)^k}{k!\,\Gamma(k+a+3/2)}
\eq
with $a=l$, and we get
\beq \label{e38}
\il_0^1\,R_n^{l,\alpha}(\rho)\,j_l(q\rho)\,\rho^2\,d\rho=\tfrac12\,\sqrt{\pi}\,(\tfrac12 q)^l \,\sum_{k=0}^{\infty}\,\frac{({-}\tfrac14 q^2)^k}{k!\,\Gamma(k+l+3/2)}\,J_{nk}^{l,\alpha}~,
\eq
where
\beq \label{e39}
J_{nk}^{l,\alpha}=\il_0^1\,\rho^{l+2k+2}\,R_n^{l,\alpha}(\rho)\,d\rho~.
\eq
We evaluate $J$ in (\ref{e39}) by using the definition of $R$ in (\ref{e5}), the substitution $t=2\rho^2-1$, Rodriguez' formula \cite{ref16}, p.~161
\beq \label{e40}
(1-t)^{\alpha}\,(1+t)^{\beta}\,P_p^{(\alpha,\beta)}(t)=\frac{({-}1)^p}{2^p\,p!}\, \Bigl(\frac{d}{dt}\Bigr)^p\,[(1-t)^{p+\alpha}\,(1+t)^{p+\beta}]~,
\eq
and subsequently $p$ partial integration. Thus
\newpage
\noindent
\begin{eqnarray} \label{e41}
J_{nk}^{l,\alpha} & = & \il_0^1\,\rho^{2l+2k+1}(1-\rho^2)^{\alpha}\,P_p^{(\alpha,l+1/2)} (2\rho^2-1)\,\rho\,d\rho \nonumber \\[3.5mm]
& = & 2^{-l-k-\alpha-5/2}\,\il_{-1}^1\,(1+t)^k\,(1-t)^{\alpha}\,(1+t)^{l+1/2}\,P_p^{(\alpha,l+1/2)} (t)\,dt \nonumber \\[3.5mm]
& = & \frac{1}{p!}\,2^{-p-l-k-\alpha-5/2}\,\il_{-1}^1\,\Bigl(\frac{d}{dt}\Bigr)^p\,[(1+t)^k]\,(1-t)^{p+\alpha} \,(1+t)^{p+l+1/2}\,dt~. \nonumber \\
\mbox{}
\end{eqnarray}
It follows that $J$ in (\ref{e39}) vanishes for $k=0,1,...,p-1\,$, and that
\beq \label{e42}
J_{nk}^{l,\alpha}=\frac{2^{-p-l-k-\alpha-5/2}}{p!\,(k-p)!}\,k!\,\il_{-1}^1\,(1-t)^{p+\alpha}\,(1+t)^{k+l+1/2}\,dt
\eq
for $k=p,p+1,...\,$. For the remaining integral, we use
\beq \label{e43}
\il_{-1}^1\,(1-t)^{r-1}\,(1+t)^{s-1}\,dt=2^{r+s-1}\,\frac{\Gamma(r)\,\Gamma(s)}{\Gamma(r+s)}
\eq
(Beta-integral), and then
\beq \label{e44}
J_{nk}^{l,\alpha}=\frac12~\frac{k!}{(k-p)!\,p!}~\frac{\Gamma(p+\alpha+1)\,\Gamma(k+l+3/2)} {\Gamma(k+l+p+\alpha+5/2)}
\eq
for $k=p,p+1,...\,$, and $J_{nk}^{l,\alpha}=0$ otherwise. When we insert this into (\ref{e38}), using the generalized Pochhammer symbol and noting various cancellations, we get
\begin{eqnarray} \label{e45}
& \mbox{} & \il_0^1\,R_n^{l,\alpha}(\rho)\,j_l(qp)\,\rho^2\,d\rho \nonumber \\[1mm]
& & =~\tfrac14\,\sqrt{\pi}\,(\tfrac12 q)^l\,(p+1)_{\alpha}\,\sum_{k=p}^{\infty}\,\frac{({-}\tfrac14 q^2)^k} {(k-p)!\,\Gamma(k+l+p+\alpha+5/2)}~.
\end{eqnarray}
We finally shift the summation index $k$ by $p$ positions, and note that $n=l+2p$, to get
\newpage
\noindent
\begin{eqnarray} \label{e46}
& \mbox{} & \il_0^1\,R_n^{l,\alpha}(\rho)\,j_l(qp)\,\rho^2\,d\rho \nonumber \\[1mm]
& & =~\tfrac12\,({-}1)^p\,(p+1)_{\alpha}\cdot\tfrac12\,\sqrt{\pi}\,(\tfrac12 q)^n\,\sum_{k=0}^{\infty}\,\frac{({-}\tfrac14 q^2)^k}{k!\,\Gamma(k+n+\alpha+1+3/2)} \nonumber \\[3.5mm]
& & =~2^{\alpha}({-}1)^p\,(p+1)_{\alpha}\,\frac{j_{n+\alpha+1}(q)}{q^{\alpha+1}}~,
\end{eqnarray}
where (\ref{e37}) has been used with $a=n+\alpha+1$. This is (\ref{e34}).

\subsection{Fourier transform of separable functions} \label{subsec2.5}
\mbox{} \\[-9mm]

With $F=AB$ as in (\ref{e14}), and expanded as in (\ref{e15}), we have
\beq \label{e47}
{\cal F}\,[F](\bfx)=\sum_{l=0}^{\infty}~\sum_{n=l,l+2,...}~\sum_{m={-}l}^l\,c_{nl}^{m,\alpha}\, {\cal F}\,[Z_{nl}^{m,\alpha}](\bfx)~,
\eq
with $c_{nl}^{m,\alpha}$ given in (\ref{e16}--\ref{e18}). For the case that $A(\vart,\varp)=h(\vart)$, with $h$ a smooth function supported by a small interval $[0,\vart_0]$, we have
\beq \label{e48}
a_l^m(A)=\delta_{m0}\,\il_0^{\vart_0}\,h(\vart)\,P_l(\cos\vart)\sin\vart\,d\vart~,
\eq
and it should, in general, be a relatively light effort to find the latter integrals. Alternatively, in \cite{ref18}, Sec.~3, a procedure is given, using scaling theory for 2D Zernike polynomials, for expressing the integrals in (\ref{e48}) in terms of expansion coefficients of $h(2\arcsin[\rho\sin(\tfrac12\vart_0)])$ with respect to ${}^2R_{2l}^{0,0}(\rho)$, $l=0,1,...\,$. By smoothness of $h$, only relatively few of the latter coefficients need to be computed.

When $F$ in (\ref{e47}) is replaced by $A_{{\rm Funk}}B$, we just need to replace $a_l^m(A)$ by $2\pi\,P_l(0)\,a_l^m(A)$, see Subsec.~\ref{subsec2.3}.

\section{Scaling theory for generalized 3D Zernike functions} \label{sec3}
\mbox{} \\[-9mm]

We present a result on scaling the radial part of the generalized 3D Zernike functions. The proof uses the same steps as the one given for the corresponding result in 2D in \cite{ref6}, Sec.~6, and starts from the following integral representation of $R_n^{l,\alpha}$.

\subsection{Integral representation of $R_n^{l,\alpha}$} \label{subsec3.1}
\mbox{} \\[-9mm]

We have for $\alpha>{-}1$ and $n=l+2p$ with $l,p=0,1,...$
\beq \label{e49}
R_n^{l,\alpha}(\rho)=\frac{2}{\pi}\,({-}1)^p\,2^{\alpha}(p+1)_{\alpha}\,\il_0^{\infty}\,\frac{j_{n+\alpha+1} (q)\,j_l(qp)}{q^{\alpha-1}}\,dq~,~~~~~~0\leq\rho<1~.
\eq
{\bf Proof.}~~By Fourier inversion, see (\ref{e36}), we have for $\bom\in B$, $\bom=\rho\bfeta$ with $0\leq\rho<1$ and $|\bfeta|=1$,
\begin{eqnarray} \label{e50}
& \mbox{} & Z_{nl}^{m,\alpha}(\nu,\mu,\sigma)=Z_{nl}^{m,\alpha}(\bom) \nonumber \\[3mm]
& & =~\il_0^{\infty}\iint\limits_{\bxi\in S}\,e^{-2\pi i\bom\cdot r\bxi}\,4\pi i^n\,2^{\alpha}(p+1)_{\alpha}\,\frac{j_{n+\alpha+1}(2\pi r)}{(2\pi r)^{\alpha+1}}\, Y_l^m(\bxi)\,r^2\,dr\,d\bxi \nonumber \\[3.5mm]
& & =~4\pi i^n\,2^{\alpha}(p+1)_{\alpha}\,\il_0^{\infty}\,\frac{j_{n+\alpha+1}(2\pi r)}{(2\pi r)^{\alpha+1}}\,\left(~\il_{\bxi\in S}\,e^{-2\pi ir\bom\cdot\bxi}\,Y_l^m(\bxi)\,d\bxi\right)\,r^2\,dr \nonumber \\[3.5mm]
& & =~4\pi i^n\,2^{\alpha}(p+1)_{\alpha}\,\il_0^{\infty}\,\frac{j_{n+\alpha+1}(2\pi r)}{(2\pi r)^{\alpha+1}}\,4\pi i^l\,Y_l^m(\bfeta)\,j_l({-}2\pi\rho r)\,r^2\,dr~,
\end{eqnarray}
where (\ref{e32}) has been used. Substituting $q=2\pi r$, using $j_l({-}z)=({-}1)^l\,j_l(z)$ and $n=l+2p$, we then get
\beq \label{e51}
Z_l^{m,\alpha}(\rho\bfeta)=\frac{2}{\pi}\,({-}1)^p\,2^{\alpha}(p+1)_{\alpha}\,Y_l^m(\bfeta)\,\il_0^{\infty} \,\frac{j_{n+\alpha+1}(q)\,j_l(qp)}{q^{\alpha-1}}\,dq~,
\eq
and this is the required result.

\subsection{Scaling the radial part} \label{subsec3.2}
\mbox{} \\[-9mm]

We have for $\alpha>{-}1$, $n=l+2p$ with $l,p=0,1,...$ and $0\leq\eps\leq1$
\beq \label{e52}
R_n^{l,\alpha}(\eps\rho)=\sum_{n'=l,l+2,...}\,C_{nn'}^{l,\alpha}(\eps)\,R_{n'}^{l,0}(\rho)~,~~~~~~0\leq\rho<1~,
\eq
where
\begin{eqnarray} \label{e53}
C_{nn'}^{l,\alpha} & = & \frac{2}{\pi}\,({-}1)^{\frac{n+n'-2l}{2}}\,2^{\alpha}(p+1)_{\alpha} \,\left[\il_0^{\infty}\,\frac{j_{n+\alpha+1}(q)\,j_{n'}(q\eps)}{q^{\alpha-1}}\,dq\right. \nonumber \\[3.5mm]
& & \hspace*{4cm}\left.+~\il_0^{\infty}\,\frac{j_{n+\alpha+1}(q)\,j_{n'+2}(q\eps)}{q^{\alpha-1}}\,dq \right]~.
\end{eqnarray}
In the case that $\alpha=0$, we have
\beq \label{e54}
C_{nn'}^{l,0}(\eps)=\left\{\ba{llll}
R_n^{n',0}(\eps)-R_n^{n'+2,0}(\eps) & \!\!, & ~~~n'=l,l+2,...,n \\[3mm]
0 & \!\!, & ~~~n'=n+2,n+4,... & \!\!,
\ea\right.
\eq
where we set $R_n^{n+2,0}=0$ for the first option in (\ref{e54}) with $n'=n$. \\ \\
{\bf Proof.}~~We have by completeness and orthogonality
\beq \label{e55}
C_{nn'}^{l,\alpha}(\eps)=(2n'+3)\,\il_0^1\,R_n^{l,\alpha}(\eps\rho)\,R_{n'}^{l,0}(\rho)\,\rho^2\,d\rho~.
\eq
We use the results (\ref{e34}) and (\ref{e49}) to write
\begin{eqnarray} \label{e56}
& \mbox{} & \il_0^1\,R_n^{l,\alpha}(\eps\rho)\,R_{n'}^{l,0}(\rho)\,\rho^2\,d\rho \nonumber \\[3.5mm]
& & =~\frac{2}{\pi}\,({-}1)^p\,2^{\alpha}(p+1)_{\alpha}\,\il_0^1\,\left(\il_0^{\infty}\,\frac{j_{n+\alpha+1} (q)\,j_l(q\eps\rho)}{q^{\alpha-1}}\,dq\right)\,R_{n'}^{l,0}(\rho)\,\rho^2\,d\rho \nonumber \\[3.5mm]
& & =~\frac{2}{\pi}\,({-}1)^p\,2^{\alpha}(p+1)_{\alpha}\,\il_0^{\infty}\,\frac{j_{n+\alpha+1} (q)}{q^{\alpha-1}}\,\left(\il_0^1\,R_n^{l,0}(\rho)\,j_l(q\eps\rho)\,\rho^2\,d\rho\right)\,dq \nonumber \\[3.5mm]
& & =~\frac{2}{\pi}\,({-}1)^{p+\frac12 (n'-l)}\,2^{\alpha}(p+1)_{\alpha}\,\il_0^{\infty}\, \frac{j_{n+\alpha+1}(q)}{q^{\alpha-1}}~\frac{j_{n'+1}(q\eps)}{q\eps}\,dq~.
\end{eqnarray}
Next, we use $p=\frac12 (n-l)$ and
\beq \label{e57}
\frac{j_{n'+1}(z)}{z}=\frac{1}{2n'+3}\,(j_{n'}(z)+j_{n'+2}(z))~,
\eq
and this gives (\ref{e53}).

The remaining integrals in (\ref{e53}) are of the form
\beq \label{e58}
\il_0^{\infty}\,\frac{j_{n+\alpha+1}(q)\,j_{n''}(q\eps)}{q^{\alpha-1}}\,dq~,~~~~~~n''=n'~{\rm or}~n'+2=l,l+2,...~,
\eq
and can be evaluated in terms of ${}_2F_1$ using the general discontinuous Weber-Schafheitlin integral \cite{ref14}, p.~244. For the case that $\alpha=0,1,...\,$, the identity in (\ref{e52}) is a relation between a polynomial comprising the powers $\rho^l,\rho^{l+2},...,\rho^{n+2\alpha}$ at the left-hand side and orthogonal polynomials $R_{n'}^{l,0}$ at the right-hand side. Therefore, $C_{nn'}^{l,\alpha}(\eps)=0$ for $n'>n+2\alpha$. In the case that $\alpha=0$, the integrals in (\ref{e58}) can all be expressed, see (\ref{e49}), as
\beq \label{e59}
\il_0^{\infty}\,j_{n+1}(q)\,j_{n''}(q\eps)\,q\,dq=\frac{\pi}{2}\,({-}1)^{\frac{n''-n}{2}}\, R_n^{n''}(\eps)~,
\eq
except the one with $n''=n+2$. The results (\ref{e49}), (\ref{e59}) can also be obtained using the Weber-Schafheidlin integral result: there holds
\begin{eqnarray} \label{e60}
& \mbox{} & \il_0^{\infty}\,\frac{j_{n+\alpha+1}(q)\,j_{n''}(q\eps)}{q^{\alpha-1}}\,dq \nonumber \\[2mm]
& & =~\frac{\pi}{2}~\frac{\Gamma(\tfrac12(n+n''+3))\,\eps^{n''}} {2^{\alpha}\,\Gamma(n''+3/2)\,\Gamma(\tfrac12(n-n'')+\alpha+1)} \nonumber \\[3.5mm]
& & \hspace*{5mm}\cdot\:{}_2F_1({-}\tfrac12 (n-n''-\alpha),\tfrac12(n+n''+3);n''+3/2;\eps^2)~.
\end{eqnarray}
This vanishes for the case that $\alpha=0,1,...$ when $n''=n+2\alpha+2,n+2\alpha+4,...$ due to the $\Gamma(\tfrac12(n-n'')+\alpha+1)$ in the denominator. Therefore, (\ref{e59}) also holds for $n''=n+2$, when we interpret the right-hand side as 0. This then yields (\ref{e54}) upon carefully keeping track of the various $({-}1)^j$.

\section{Radon transform and recursions for computing the radial parts of generalized 3D Zernike functions} \label{sec4}
\mbox{} \\[-9mm]

We present a closed-form expression for the Radon transform
\newpage
\beq \label{e61}
({\cal R} Z)(\tau,\bfeta_R)=\il_{\bom\in B}\,Z(\bom)\,\delta(\tau-\bom\cdot\bfeta_R)\,d\bom~,
\eq
with $\tau\in\dR$ and $\bfeta_R\in S$, of any generalized 3D Zernike function $Z=Z_{nl}^{m,\alpha}$. This result is due to A.K.\ Louis, Theorem~3.1, case $N=3$, in \cite{ref7}, and can be proved using the projection theorem for Radon transforms, the Funk-Hecke theorem, the explicit form of the Fourier transform of generalized 3D Zernike functions, and an integral result for the Fourier transform of the right-hand side of (\ref{e34}). See also the proof for the 2D case in \cite{ref6}, Sec.~5. This is then used to derive an integral representation of the radial part of any generalized 3D Zernike function, which can be employed to find a recursion of the Shakibaei-Paramesran type \cite{ref13} for computing all $R_n^{l,\alpha}(\rho)$, $l,n=0,1,...$ and $n-l$ even and non-negative, at a particular $\rho\in[0,1]$ from $R_n^{l,\alpha}(\rho)$ with $l=n=0$.

\subsection{Radon transform of generalized 3D Zernike functions} \label{subsec4.1}
\mbox{} \\[-9mm]

We have for integer $n,l,m$ with $n-l$ even and non-negative and $n,l\geq0$, and $m={-}l,{-}l+1,...,l\,$, all $\alpha>{-}1$ and all $\tau\in\dR$, $\bfeta_R\in S$
\beq \label{e62}
({\cal R}Z_{nl}^{m,\alpha})(\tau,\bfeta_R)=\frac1c\,(1-\tau^2)^{\alpha+1}\,C_n^{\alpha+3/2}(\tau)\,Y_l^m(\bfeta_R)~,
\eq
where $C_n^{\alpha+3/2}$ is the Gegenbauer polynomial, \cite{ref14}--\cite{ref16}, of degree $n$ corresponding to the weight $(1-\tau^2)^{\alpha+1}$, $|\tau|\leq1$, and
\beq \label{e63}
c=\frac{2^{-2(1+\alpha)}}{\sqrt{\pi}\,(p+1)_{\alpha}}~\frac{\Gamma(n+2\alpha+3)}{\Gamma(n+1)\, \Gamma(\alpha+3/2)}~.
\eq
This follows from \cite{ref7}, Theorem~3.1 (choice $N=3$, $\nu=\alpha+3/2$, $w_{\nu}(s)=(1-s^2)^{\alpha+1}$ in terms of the parameters and weights in \cite{ref7}).

\subsection{Integral representation of $R_n^{l,\alpha}$} \label{subsec4.2}
\mbox{} \\[-9mm]

We have for any integer $n,l\geq0$ such that $n-l$ is even and non-negative, any $\alpha>{-}1$, and $0\leq\rho\leq1$
\beq \label{e64}
R_n^{l,\alpha}(\rho)=\dfrac12~\frac{(3/2)_{p+l}}{(\alpha+3/2)_{p+l}}\,(1-\rho^2)^{\alpha}\,\il_{-1}^1\, C_n^{\alpha+3/2}(\rho t)\,P_l(t)\,dt~.
\eq
{\bf Proof.}~~We develop the function
\beq \label{e65}
\bom=(\nu,\mu,\sigma)\in B\mapsto(1-\rho^2)^{\alpha}\,C_n^{\alpha+3/2}(\sigma)
\eq
into a $Z^{\alpha}$-series, so that
\beq \label{e66}
(1-\rho^2)^{\alpha}\,C_n^{\alpha+3/2}(\sigma)=\sum_{n',l,m}\,\beta_{n'l}^m\,Z_{n'l}^{m,\alpha}(\nu,\mu,\sigma)~,
\eq
with
\begin{eqnarray} \label{e67}
\beta_{n'l}^m & = &
\frac{1}{N_{n'l}^{\alpha}}~\:\iiint\limits_{\nu^2+\mu^2+\sigma^2\leq1}\,C_n^{\alpha+3/2}(\sigma) (Z_{n'l}^{m,\alpha}(\nu,\mu,\sigma))^{\ast}\,d\nu d\mu d\sigma \nonumber \\[3.5mm]
& = & \frac{1}{N_{n'l}^{\alpha}}\,\il_{-1}^1\,C_n^{\alpha+3/2}(\sigma)\left(~\:\iint\limits _{\nu^2+\mu^2\leq1-\sigma^2}\,Z_{n'l}^m(\nu,\mu,\sigma)\,d\nu d\mu\right)^{\ast}\,d\sigma~.
\end{eqnarray}
To evaluate the double integral in the last member of (\ref{e67}), we use (\ref{e62}) with $\bfeta_R=(0,0,1)$ and $\tau=\sigma$, so that
\begin{eqnarray} \label{e68}
\iint\limits_{\nu^2+\mu^2\leq1-\sigma^2}\,Z_{n'l}^{m,\alpha}(\nu,\mu,\sigma)\,d\nu d\mu & = & \il_{\bom\in B}\,Z_{n'l}^{m,\alpha}(\bom)\,\delta_0(\tau-\bom\cdot\bfeta_R)\,d\bom \nonumber \\[3.5mm]
& = & \frac1c\,(1-\sigma^2)^{\alpha+1}\,C_{n'}^{\alpha+3/2}(\sigma)\,Y_l^m((0,0,1))~. \nonumber \\
\mbox{}
\end{eqnarray}
Now $\bom=(0,0,1)$ when $\vart=0$, $\varp$ arbitrary, $\rho=1$ in (\ref{e4}), and so
\beq \label{e69}
Y_l^m((0,0,1))=\Bigl(\frac{2l+1}{4\pi}\Bigr)^{1/2}\,P_l^0(1)\,\delta_{m0}=\Bigl(\frac{2l+1}{4\pi}\Bigr)^{1/2}\,\delta_{0m}
\eq
by (\ref{e6}). Hence
\beq \label{e70}
\iint\limits_{\nu^2+\mu^2\leq1-\sigma^2}\,Z_{n'l}^{m,\alpha}(\nu,\mu,\sigma)\,d\nu d\mu=\frac1c\,(1-\sigma^2)^{\alpha+1}\,C_{n'}^{\alpha+3/2}(\sigma)\Bigl(\frac{2l+1}{4\pi}\Bigr)^{1/2}\,\delta_{m0}~.
\eq
This yields, by orthogonality of the $C^{\alpha+3/2}$,
\begin{eqnarray} \label{e71}
\beta_{n'l}^m & = & \frac{\delta_{m0}}{c\,N_{n'l}^{\alpha}}\,\Bigl(\frac{2l+1}{4\pi}\Bigr)^{1/2}\,\il_{-1}^1\,C_n^{\alpha+3/2}(\sigma)\,C_{n'}^{\alpha+3/2}(\sigma)(1-\sigma^2)^{\alpha+1} \,d\sigma \nonumber \\[3.5mm]
& = & \frac{\delta_{m0}\,\delta_{nn'}}{c\,N_{nl}^{\alpha}}\,\Bigl(\frac{2l+1}{4\pi}\Bigr)^{1/2}\,M_n^{\alpha}~,
\end{eqnarray}
where
\beq \label{e72}
M_n^{\alpha}=\il_{-1}^1\,(1-\sigma^2)^{\alpha+1}\,|C_n^{\alpha+3/2}(\sigma)|^2\,d\sigma= \frac{\pi\,2^{-2\alpha-2}\,\Gamma(n+2\alpha+3)}{n!\, (n+\alpha+3/2)\,\Gamma^2(\alpha+3/2)}
\eq
by \cite{ref16}, (7.8) on p.~279 with $\lambda=\alpha+3/2$. Remembering the definitions of $N_{nl}^{\alpha}$ in (\ref{e13}) and of $c$ in (\ref{e63}), we get
\beq \label{e73}
\beta_{n'l}^m=\pi(2l+1)^{1/2}\,\frac{(p+l+3/2)_{\alpha}}{\Gamma(\alpha+3/2)}\,\delta_{m0}\,\delta_{nn'}~.
\eq
It follows, see (\ref{e3}--\ref{e6}), with $\sigma=\rho\cos\vart$,
\begin{eqnarray} \label{e74}
& \mbox{} & (1-\rho^2)^{\alpha}\,C_n^{\alpha+3/2}(\rho\cos\vart) \nonumber \\[3.5mm]
& & =~\frac{\pi}{\Gamma(\alpha+3/2)}\,\sum_l\,(2l+1)^{1/2}\,(p+l+3/2)_{\alpha}\,Z_{nl}^{0,\alpha}(\nu,\mu,\sigma) \nonumber \\[3.5mm]
& & =~\frac{\sqrt{\pi}}{2\Gamma(\alpha+3/2)}\,\sum_l\,(2l+1)(p+l+3/2)_{\alpha}\,R_n^{l,\alpha}(\rho)\,P_l(\cos\vart)~.
\end{eqnarray}
From
\beq \label{e75}
\il_0^{\pi}\,P_{l_1}(\cos\vart)\,P_{l_2}(\cos\vart)\sin\vart\,d\vart=\frac{\delta_{l_1l_2}}{l_{1,2}+1/2}~,
\eq
we then get for $l=n,n-2,...,n-2\lfloor n/2\rfloor$
\begin{eqnarray} \label{e76}
& \mbox{} & \il_0^{\pi}\,(1-\rho^2)^{\alpha}\,C_n^{\alpha+3/2}(\rho\cos\vart)\,P_l(\cos\vart)\sin\vart\,d\vart \nonumber \\[3.5mm]
& & =~\frac{\sqrt{\pi}}{\Gamma(\alpha+3/2)}\,(p+l+3/2)_{\alpha} \,R_n^{l,\alpha}(\rho)~.
\end{eqnarray}
Finally, using
\beq \label{e77}
\frac{1}{\sqrt{\pi}}~\frac{\Gamma(\alpha+3/2)}{(p+l+3/2)_{\alpha}}=\frac12~\frac{(3/2)_{p+l}}{(\alpha+3/2)_{p+l}}~,
\eq
and substituting $t=\cos\vart\in[{-}1,1]$ in the integral in (\ref{e76}), we get (\ref{e64}). \\ \\
{\bf Note.}~~In \cite{ref6}, Theorem~5.2, the 2D version of this result has been given, and this leads to an expression for ${}^2R_n^{|m|,\alpha}(\rho)$ in terms of the Fourier coefficients of $C_n^{\alpha+1}(\rho\cos\vart)$. By discretization $\vart=\frac{2\pi k}{N}$ with $N$ sufficiently large, this yields a DFT-scheme for computing all $R_n^{|m|}(\rho)$, integer $m$ such that $n-|m|$ is even and non-negative, for a fixed $n=0,1,...$ and $\rho\in[0,1]$. Such a thing is more awkward in the 3D because of the occurrence of a Legendre polynomial, rather than a Chebyshev polynomial in (\ref{e64}).

\subsection{Recursion for $R_n^{l,\alpha}(\rho)$} \label{subsec4.3}
\mbox{} \\[-9mm]

We let
\beq \label{e78}
I_{nl}^{\alpha}=I_{nl}^{\alpha}(\rho)=\tfrac12\,\il_{-1}^1\,C_n^{\alpha+3/2}(\rho t)\,P_l(t)\,dt~,
\eq
so that
\beq \label{e79}
\rho^l\,P_p^{(\alpha,l+1/2)}(2\rho^2-1)=(1-\rho^2)^{-\alpha}\,R_n^{l,\alpha}(\rho)=\frac{(3/2)_{p+l}}{(\alpha+3/2)_{p+l}}\,I_{nl}^{\alpha}(\rho)~.
\eq
Then for integer $n,l\geq0$ with $n+1-l$ even and non-negative,
\beq \label{e80}
I_{n+1,l}^{\alpha}=\frac{n+\alpha+3/2}{n+1}\,\rho\,\Bigl[\frac{l+1}{l+1/2}\,I_{n,l+1}^{\alpha}+\frac{l}{l+1/2}\,I_{n,l-1}^{\alpha}\Bigr]- \frac{n+2\alpha+2}{n+1}\,I_{n-1,l}^{\alpha}~.
\eq
With the initialization $I_{00}^{\alpha}(\rho)=1$, $I_{nl}^{\alpha}\equiv0$ when $n<l$, all $I_{nl}^{\alpha}(\rho)$ can be computed for a fixed $\rho\in[0,1]$ according to the scheme pictured in Fig.~1, where the large-size numbers $0,1,2,...$ indicate the order of computation. A similar recursion exists in the 2D case. Setting
\beq \label{e81}
J_{nm}^{\alpha}=J_{nm}^{\alpha}(\rho)=\frac{(\alpha+1)_{p+m}}{(1)_{p+m}}\,(1-\rho^2)^{-\alpha}\,{}^2R_n^{m,\alpha}(\rho)
\eq
for integer $n,m\geq0$ such that $n-m$ is even and non-negative, the following holds. Let $n$, $m$ be non-negative integers such that $n+1-m$ is even and non-negative, and let $\rho\in[0,1]$. Then
\beq \label{e82}
J_{n+1,m}^{\alpha}=\frac{n+\alpha+1}{n+1}\,[J_{n,|m-1|}^{\alpha}+J_{n,m+1}^{\alpha}]-\frac{n+2\alpha+1}{n+1}\,J_{n-1,m}^{\alpha}~,
\eq
where we initialize with $J_{00}^{\alpha}(\rho)=1$, $J_{nm}^{\alpha}\equiv0$ when $n<m$. This generalizes the recursive scheme
\beq \label{e83}
{}^2R_n^{|m|}(\rho)=\rho\,[{}^2R_{n-1}^{|m-1|}(\rho)+{}^2R_{n-1}^{|m+1|}(\rho)]-{}^2R_{n-2}^{|m|}(\rho)~,
\eq
valid for the radial parts of the standard 2D Zernike circle polynomials \cite{ref13}.

\begin{figure} [!h]
\centerline{\includegraphics[width=0.6\textwidth]{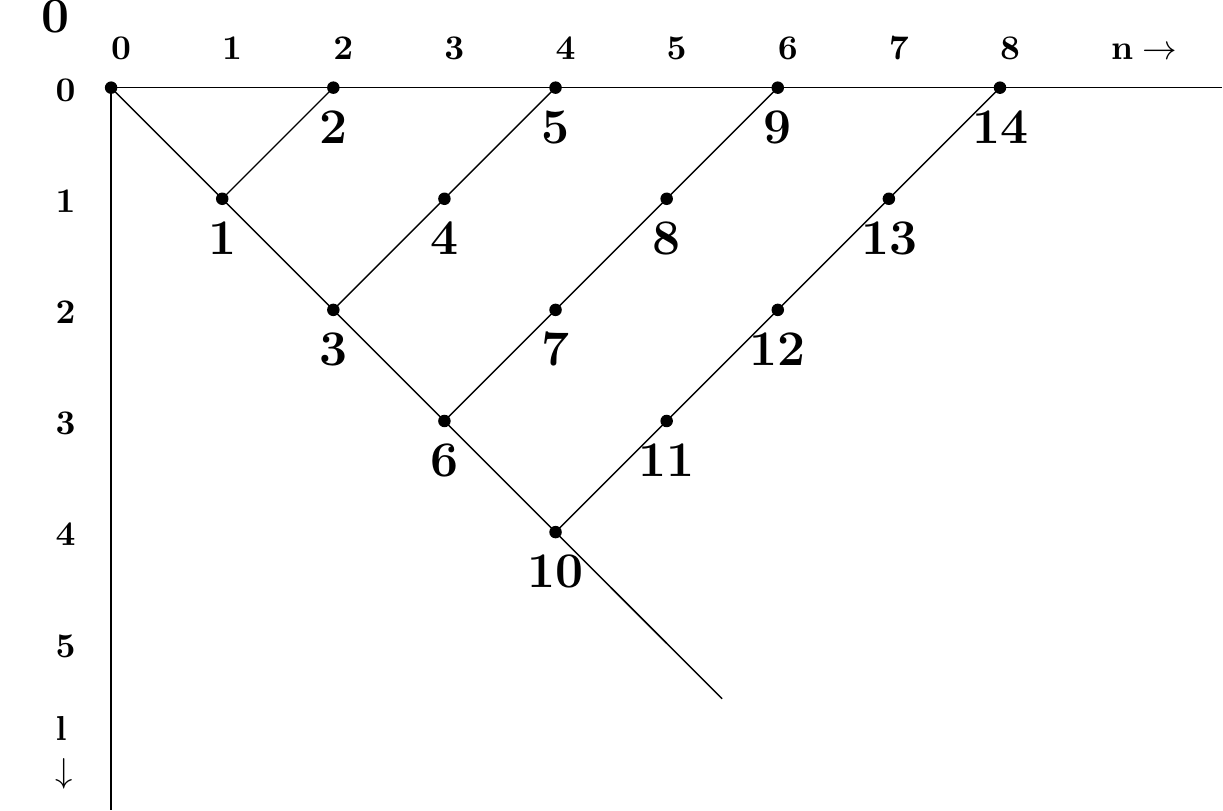}}
\caption{Order of computation (large-size numbers) of all $I_{nl}^{\alpha}(\rho)$, starting from $I_{00}^{\alpha}(\rho)=1$ and $I_{nl}^{\alpha}(\rho)=0$ for $n<l$, for a fixed $\rho\in[0,1]$ using the recursion (\ref{e80}).}
\end{figure}

To show (\ref{e80}), we use the recursions \cite{ref16}, (7.7) on p.~279 and (10.2) on p.~190,
\beq \label{e84}
(n+1)\,C_{n+1}^{\lambda}(t)=2(n+\lambda)\,t\,C_n^{\lambda}(t)-(n+2\lambda-1)\,C_{n-1}^{\lambda}(t)~,
\eq
\beq \label{e85}
(2n+1)\,t\,P_n(t)=(n+1)\,P_{n+1}(t)+n\,P_{n-1}(t)~,
\eq
valid for $n=0,1,...\,$, where we set $C_{-1}^{\lambda}(t)=0=P_{-1}(t)$. Thus from (\ref{e78}), with $n+1$ instead of $n$ and $\lambda=\alpha+3/2$,
\begin{eqnarray} \label{e86}
& \mbox{} & I_{n+1,l}^{\alpha}\nonumber \\[2mm]
& & =~\tfrac12\,\il_{-1}^1\,C_{n+1}^{\alpha+3/2}(\rho t)\,P_l(t)\,dt \nonumber \\[3.5mm]
& & =~\tfrac12\,\il_{-1}^1\,\Bigl(\frac{2(n+\alpha+3/2)}{n+1}\,\rho t\,C_n^{\alpha+3/2}(\rho t)- \frac{n+2\alpha+2}{n+1}\,C_{n-1}^{\alpha+3/2}(\rho t)\Bigr)\, P_l(t)\,dt \nonumber \\[3.5mm]
& & =~\frac{n+\alpha+3/2}{n+1}\,\rho\,\il_{-1}^1\,C_n^{\alpha+3/2}(\rho t)\,t\,P_l(t)\,dt-\frac{n+2\alpha+2}{n+1}\,I_{n-1,l}^{\alpha} \nonumber \\[3.5mm]
& & =~\frac{n+\alpha+3/2}{n+1}\,\rho\,\il_{-1}^1\,C_n^{\alpha+3/2}(\rho t)\,\Bigl(\frac{l+1}{2l+1}\,P_{l+1}(t)+\frac{l}{2l+1}\,P_{l-1}(t)\Bigr)\,dt \nonumber \\[3.5mm]
& & \hspace*{5mm}-~\frac{n+2\alpha+2}{n+1}\,I_{n-1,l}^{\alpha}~,
\end{eqnarray}
and this yields the right-hand side of (\ref{e80}) via (\ref{e78}) with $l-1$ instead of $l$.

\section{Recursion for the radial part of expansion coefficients of separable functions on the unit ball} \label{sec5}
\mbox{} \\[-9mm]

We present a recursion relation between the radial part coefficients $b(B)$ of a separable function $F=AB$ on the unit ball, see Subsec.~\ref{subsec2.2}. We recall that for integer $n,l\geq0$ such that $n-l$ is even and non-negative
\beq \label{e87}
b_n^{l,\alpha}(B)=\il_0^1\,B(\rho)\,R_n^{l,\alpha}(\rho)\,\frac{\rho^2\,d\rho}{(1-\rho^2)^{\alpha}}~.
\eq
We shall explicitly find the connection coefficients $C$ in the formula
\beq \label{e88}
R_{l+2+2p}^{l+2,\alpha}=\sum_{s=0}^{p+1}\,C_{ps}^{l+2\pr l,\alpha}\,R_{l+2s}^{l,\alpha}~,
\eq
where $l,p=0,1,...\,$. As a consequence, we have
\beq \label{e89}
b_{l+2+2p}^{l+2,\alpha}(B)=\sum_{s=0}^{p+1}\,C_{ps}^{l+2\pr l,\alpha}\,b_{l+2s}^{l,\alpha}(B)~.
\eq
Hence, from $b_n^{0,\alpha}(B)$ and $b_n^{1,\alpha}(B)$, all required $b$-coefficients can be computed recursively using (\ref{e89}).

To find the $C$'s in (\ref{e88}), we use the definition of $R$ in (\ref{e5}) and set $t=2\rho^2-1$, so that (\ref{e68}) becomes
\beq \label{e90}
\tfrac12\,(1+t)\,P_p^{(\alpha,l+5/2)}(t)=\sum_{s=0}^{p+1}\,C_{ps}^{l+2\pr l,\alpha}\,P_s^{(\alpha,l+1/2)}(t)~.
\eq
We have, see \cite{ref16}, 2$^{{\rm nd}}$ item in (4.16) on p.~166,
\beq \label{e91}
\tfrac12\,(1+t)(2n+\alpha+\beta+2)\,P_n^{(\alpha,\beta+1)}(t)=(n+\beta+1)\,P_n^{(\alpha,\beta)}(t)+(n+1)\,P_{n+1}^{(\alpha,\beta)}(t)~,
\eq
and when this is used with $p=n$, $\alpha=\alpha$, $\beta=l+1/2$, we get
\begin{eqnarray} \label{e92}
& \mbox{} &
\hspace*{-8mm}\tfrac12\,(1+t)\,P_p^{(\alpha,l+5/2)}(t) \nonumber \\[3.5mm]
& & \hspace*{-8mm}=~\frac{p+l+5/2}{2p+\alpha+l+7/2}\,P_p^{(\alpha,l+3/2)}(t)+\frac{p+1}{2p+\alpha+l+7/2}\,P_{p+1}^{(\alpha,l+3/2)}(t)~.
\end{eqnarray}
Next, we have, see \cite{ref19}, (7.32),
\begin{eqnarray} \label{e93}
& \mbox{} & \hspace*{-5mm}P_n^{(\alpha,\delta)}=\frac{(\alpha+1)_n}{(\alpha+\beta+2)_n} \nonumber \\[3.5mm]
& & \hspace*{-5mm}\cdot\:\sum_{k=0}^n\,\frac{({-}1)^{n-k}(\delta-\beta)_{n-k}(\alpha+\beta+1)_k(\alpha+\beta+2)_{2k}(n+\alpha+\delta+1)_k} {(1)_{n-k}(\alpha+1)_k(\alpha+\beta+1)_{2k}(n+\alpha+\beta+2)_k}\,P_k^{(\alpha,\beta)}~, \nonumber \\
\mbox{}
\end{eqnarray}
and when we use this with $n=p,p+1$, $\alpha=\alpha$, $\beta=l+1/2$, $\delta=\beta+1=l+3/2$, we get
\begin{eqnarray} \label{e94}
& \mbox{} & P_n^{(\alpha,l+3/2)} \nonumber \\[2mm]
& & =~\frac{(\alpha+1)_n}{(\alpha+l+5/2)_n}\,\sum_{k=0}^n\,({-}1)^{n-k}\,\frac{(\alpha+l+3/2)_k(\alpha+l+5/2)_{2k}} {(\alpha+1)_k(\alpha+l+3/2)_{2k}}\,P_k^{(\alpha,l+1/2)}~. \nonumber \\
\mbox{}
\end{eqnarray}
From (\ref{e92}) and (\ref{e94}), we compute then
\begin{eqnarray} \label{e95}
& \mbox{} & \hspace*{-8mm}C_{ps}^{l+2\pr l,\alpha} \nonumber \\[2mm]
& & \hspace*{-8mm}=~\Bigl(\frac{p+l+5/2}{2p+\alpha+l+7/2}~\frac{(\alpha+1)_p}{(\alpha+l+5/2)_p}-\frac{p+1}{2p+\alpha+l+7/2}~\frac{(\alpha+1)_{p+1}}{(\alpha+l+5/2)_{p+1}}\Bigr) \nonumber \\[3.5mm]
& & \hspace*{-3mm}\cdot\:({-}1)^{p-s}\,\frac{(\alpha+l+3/2)_s(\alpha+l+5/2)_{2s}}{(\alpha+1)_s(\alpha+l+3/2)_{2s}}~,~~~~~~s=0,1,...,p~,
\end{eqnarray}
and, for $s=p+1$, we get
\begin{eqnarray} \label{e96}
& \mbox{} & \hspace*{-6mm}C_{p,p+1}^{l+2\pr l,\alpha} \nonumber \\[2mm]
& & \hspace*{-6mm}=~\frac{p+1}{2p+\alpha+l+7/2}~\frac{(\alpha+1)_{p+1}}{(\alpha+l+5/2)_{p+1}}~ \frac{(\alpha+l+3/2)_{p+1}(\alpha+l+5/2)_{2p+2}}{(\alpha+1)_{p+1}(\alpha+l+3/2)_{2p+2}}~. \nonumber \\[2mm]
\mbox{}
\end{eqnarray}
For $s=0,1,...,p\,$, we further compute the factor on the first line of (\ref{e95}) as
\begin{eqnarray} \label{e97}
& \mbox{} & \frac{(\alpha+1)_p}{(2p+\alpha+l+7/2)(\alpha+l+5/2)_{p+1}} \nonumber \\[3.5mm]
& & \hspace*{2cm}\cdot\:[(p+l+5/2)(p+l+5/2+\alpha)-(p+1)(p+1+\alpha)] \nonumber \\[3.5mm]
& & =~\frac{(\alpha+1)_p}{(\alpha+l+5/2)_{p+1}}\,(l+3/2)~.
\end{eqnarray}
Hence, for $s=0,1,...,p$
\beq \label{e98}
C_{ps}^{l+2\pr l,\alpha}=({-}1)^{p-s}(l+3/2)\,\frac{(\alpha+1)_p(\alpha+l+3/2)_s(\alpha+l+5/2)_{2s}} {(\alpha+l+5/2)_{p+1}(\alpha+1)_s(\alpha+l+3/2)_{2s}}~,
\eq
and (\ref{e96}) for $s=p+1$ simplifies to
\beq \label{e99}
C_{p,p+1}^{l+2\pr l,\alpha}=\frac{p+1}{\alpha+l+p+5/2}~.
\eq
Also, see (\ref{e103}) below for simplification of the ratio of the two Pochhammer symbols of order $2s$ in (\ref{e98}).

The corresponding result for the 2D case to write ${}^2R_{l+2+2p}^{l+2,\alpha}$ as a linear combination of ${}^2R_{l+2s}^{l,\alpha}$ reads
\beq \label{e100}
{}^2R_{l+2+2p}^{l+2,\alpha}=\sum_{s=0}^{p+1}\,{}^2C_{ps}^{l+2\pr l,\alpha}\,{}^2R_{l+2s}^{l,\alpha}~,
\eq
with
\beq \label{e101}
{}^2C_{ps}^{l+2\pr l,\alpha}=({-}1)^{p-s}(l+1)\,\frac{(\alpha+1)_p(\alpha+l+1)_s}{(\alpha+l+2)_{p+1}(\alpha+1)_s}~\frac{\alpha+l+2s+1} {\alpha+l+1}
\eq
for $s=0,1,...,p$ and
\beq \label{e102}
{}^2C_{p,p+1}^{l+2\pr l,\alpha}=\frac{p+1}{\alpha+l+p+2}~.
\eq
For the case that $\alpha=0$, this result was already established in 1942 by Nijboer, see \cite{ref20}, (2.24) and \cite{ref21}, (200).

The expression (\ref{e98}) for $C$ contains Pochhammer symbols of potentially high order that can, therefore, become very large. The quantities $C$ themselves are of the order of unity. They can be computed more reliably as follows. We first write (\ref{e98}) as
\beq \label{e103}
C_{ps}^{l+2\pr l,\alpha}=({-}1)^{p-s}\,\frac{(l+3/2)(\alpha+l+2s+3/2)}{\alpha+l+3/2}\,K_{ps}~,
\eq
where for a fixed $p=0,1,...$
\beq \label{e104}
K_{ps}=\frac{(\alpha+1)_p}{(\alpha+l+5/2)_{p+1}}~\frac{(\alpha+l+3/2)_s}{(\alpha+1)_s}~,~~~~~s=0,1,...,p~.
\eq
The quantity in (\ref{e103}) in front of $K_{ps}$ is innocent. As to $K_{ps}$ itself, we observe that
\beq \label{e105}
K_{pp}=\frac{\alpha+l+3/2}{(\alpha+l+p+3/2)(\alpha+l+p+5/2)}~,
\eq
and
\beq \label{e106}
K_{p,s-1}=\frac{\alpha+s}{\alpha+l+s+1/2}\,K_{ps}~,~~~~~~s=p,p-1,...,1~.
\eq

\section{Development of scaled-and-truncated radial profiles} \label{sec6}
\mbox{} \\[-9mm]

We present a procedure to find for a given even integer $\beta\geq0$ and a $\delta>{-}1$ the expansion of
\beq \label{e107}
B(\rho;\eps)=\rho^{\beta}\,C(\eps\rho)(1-\rho^2)^{\delta}~,~~~~~~0\leq\rho\leq1~,
\eq
into a series comprising the radial functions $R_{2k}^{0,\delta}$, $k=0,1,...\,$. Here $C(\rho)$, $0\leq\rho\leq1$, is an integrable function of which we assume that there is available the expansion
\beq \label{e108}
C(\rho)=\sum_{l=0}^{\infty}\,c_l\,R_{2l}^{0,0}(\rho)~,~~~~~~0\leq\rho\leq1~,
\eq
with $R_{2l}^{0,0}$ the standard 3D Zernike polynomials, and $\eps$ is a scaling parameter with $0\leq\eps\leq1$. From the expansion of $B(\rho;\eps)$ as an $R_{2k}^{0,\delta}$-series, we can obtain, by the procedure of Sec.~\ref{sec5}, all expansions of $B(\rho;\eps)$ as an $R_{2k}^{l,\delta}$-series with $l=2,4,...\,$.

From the scaling result proved in Sec.~\ref{sec3}, we have
\beq \label{e109}
R_{2l}^{0,0}(\eps\rho)=\sum_{k=0}^l\,(R_{2l}^{2k,0}(\eps)-R_{2l}^{2k+2,0}(\eps))\,R_{2k}^{0,0}(\rho)~,
\eq
and so we find
\begin{eqnarray} \label{e110}
C(\eps\rho) & = & \sum_{l=0}^{\infty}\,c_l\,R_{2l}^{0,0}(\eps\rho) \nonumber \\[3.5mm]
& = & \sum_{l=0}^{\infty}\,c_l\,\sum_{k=0}^l\,(R_{2l}^{2k,0}(\eps)-R_{2l}^{2k+2,0}(\eps))\,R_{2k}^{0,0}(\rho) \nonumber \\[3.5mm]
& = & \sum_{k=0}^{\infty}\,c_k(\eps)\,R_{2k}^{0,0}(\rho)~,
\end{eqnarray}
where
\beq \label{e111}
c_k(\eps)=\sum_{l=k}^{\infty}\,(R_{2l}^{2k,0}(\eps)-R_{2l}^{2k+2,0}(\eps))\,c_l~.
\eq
Compare \cite{ref22}, \cite{ref23}. In particular, we have $c_k(1)=c_k$ since $R_{2i}^{2j,0}(1)=1$ or 0 according as $i\geq j$ or $i<j$. Also, $c_k(0)=\delta_{k0}\,c_0$ since $R_{2k}^{0,0}(0)=\delta_{k0}$. Furthermore, when $c_l=0$ for $l\geq L$, we have $c_k(\eps)=0$ for $k\geq L$.

The required expansion of $B(\rho;\eps)$ in (\ref{e107}) can now be obtained from the expansion of $C(\eps\rho)$ in (\ref{e110}) by expanding systematically
\beq \label{e112}
\rho^{\beta}\,R_{2k}^{0,0}(\rho)(1-\rho^2)^{\delta}=\sum_i\,D_{ki}^{\beta\delta}\,R_{2i}^{0,\delta}(\rho)~,
\eq
where we recall that $\beta=2r$ is even. We achieve this using two steps, viz.\ by expanding
\beq \label{e113}
\rho^{2r}\,R_{2k}^{0,0}(\rho)=\sum_j\,E_{kj}^r\,R_{2j}^{0,0}(\rho)~,
\eq
\beq \label{e114}
R_{2j}^{0,0}(\rho)(1-\rho^2)^{\delta}=\sum_i\,F_{ji}^{\delta}\,R_{2i}^{0,\delta}(\rho)~,
\eq
respectively.

\subsection{Computation of the expansion coefficients in (\ref{e113})} \label{subsec6.1}
\mbox{} \\[-9mm]

We have by orthogonality and (\ref{e13}), case $\alpha=0$,
\beq \label{e115}
E_{kj}^r=2(2j+3/2)\,\il_0^1\,\rho^{2r}\,R_{2k}^{0,0}(\rho)\,R_{2j}^{0,0}(\rho)\,\rho^2\,d\rho~.
\eq
We next expand $\rho^r\,R_{2k}^{0,0}(\rho)$ and $\rho^r\,R_{2j}^{0,0}(\rho)$ into a $R_{r+2l}^{r,0}$-series according to
\beq \label{e116}
\rho^r\,R_{2i}^{0,0}(\rho)=\sum_l\,C_{il}^r\,R_{r+2l}^{r,0}(\rho)~,~~~~~~i=k,j~.
\eq
Then by orthogonality and (\ref{e10}), case $\alpha=0$, we get
\beq \label{e117}
E_{kj}^r=\sum_l\,\frac{2j+3/2}{r+2l+3/2}\,C_{kl}^r\,C_{jl}^r~.
\eq
The $C_{il}^r$ can be evaluated explicitly as
\beq \label{e118}
C_{il}^r=\frac{r+2l+3/2}{i+r+l+3/2}\,\Bigl(\!\ba{c} r \\[1mm] i-l \ea\!\Bigr)\,\Bigl(\!\ba{c} i+l+1/2 \\[1mm] l \ea\!\Bigr)
\,/\,
\Bigl(\!\ba{c} i+r+l+1/2 \\[1mm] i \ea\!\Bigr)
\eq
when $l=0,1,...,i$ and $i-l\leq r$, and 0 otherwise. Thus, the summation range in (\ref{e117}) consists of all integer $l$ such that
\beq \label{e119}
\max\,\{0,k-r,j-r\}\leq l\leq\min\,\{k,j\}~.
\eq
This summation range is empty, and so $C_{kl}^r=0$, if and only if $|k-j|>r$.

To show (\ref{e118}), we can use \cite{ref19}, (7.32), compare (\ref{e93}), with $\alpha=0$, $\delta=1/2$, $\beta=r+1/2$, since (\ref{e116}) is the same as
\beq \label{e120}
P_i^{(0,1/2)}=\sum_l\,C_{il}^r\,P_l^{(0,r+1/2)}~.
\eq
In terms of Pochhammer symbols, we have
\beq \label{e121}
C_{il}^r=\Bigl(\!\ba{c} r \\[1mm] i-l \ea\!\Bigr)\,\frac{r+2l+3/2}{r+3/2}~\frac{i!\,(r+3/2)_l(i+3/2)_l}{l!\,(r+5/2)_i(i+r+5/2)_l}~,
\eq
and this can be case into the binomial form (\ref{e118}).

It does not seem possible to find the $E_{kj}^r$ in (\ref{e113}) in closed form. This problem also occurs in the 2D case. In the 2D case, one expands
\beq \label{e122}
\rho^{2r}\,{}^2R_{2k}^{0,0}(\rho)=\sum_j\,{}^2E_{kj}^r\,{}^2R_{2j}^{0,0}(\rho)
\eq
with
\beq \label{e123}
{}^2E_{kj}^r=\sum_l\,\frac{2j+1}{r+2l+1}\,{}^2C_{kl}^r\,{}^2C_{jl}^r~,
\eq
and, for $i=k,j$ with $l=0,1,...,i$ and $i-l\leq r$,
\beq \label{e124}
{}^2C_{il}^r=\frac{r+2l+1}{i+r+l+1}\,\Bigl(\!\ba{c} r \\[1mm] i-l \ea\!\Bigr)\,\Bigl(\!\ba{c} i+l \\[1mm] l \ea\!\Bigr)
\,/\,
\Bigl(\!\ba{c} r+i+l \\[1mm] i \ea\!\Bigr)~.
\eq

\subsection{Computation of the expansion coefficients in (\ref{e114})} \label{subsec6.2}
\mbox{} \\[-9mm]

We have that (\ref{e114}) is the same as
\beq \label{e125}
P_j^{(0,1/2)}=\sum_i\,F_{ji}^{\delta}\,P_i^{(\delta,1/2)}~.
\eq
We now use \cite{ref19}, 7.33,
\begin{eqnarray} \label{e126}
P_n^{(\gamma,\beta)} & \!\!= & \!\!\frac{(\beta+1)_n}{(\alpha+\beta+2)_n} \nonumber \\[3.5mm]
& & \!\!\cdot\:\sum_{k=0}^n\,\frac{(\gamma-\alpha)_{n-k}(\alpha+\beta+1)_k(\alpha+\beta+2)_{2k}(n+\gamma+\beta+1)_k} {(1)_{n-k}(\beta+1)_k(\alpha+\beta+1)_{2k}(n+\alpha+\beta+2)_k}\,P_k^{(\alpha,\beta)} \nonumber \\[2mm]
\mbox{}
\end{eqnarray}
with $\gamma=0$, $\beta=1/2$, $n=j$, $\alpha=\delta$, $k=i$. We then find
\begin{eqnarray} \label{e127}
F_{ji}^{\delta} & \!\!= & \!\!({-}1)^{j-i}\Bigl(\!\ba{c} \delta \\[1mm] j-i \ea\!\Bigr)\, \frac{\delta+2i+3/2}{\delta+3/2}~\frac{(3/2)_j}{(\delta+5/2)_j}~ \frac{(\delta+3/2)_i(j+3/2)_i}{(3/2)_i(j+\delta+5/2)_i} \nonumber \\[3.5mm]
& \!\!= & \!\!({-}1)^{j-i}\,\frac{\delta+2i+3/2}{\delta+i+j+3/2}\,\Bigl(\!\ba{c} \delta \\[1mm] j-i \ea\!\Bigr) \Bigl(\!\ba{c} i+j+1/2 \\[1mm] j \ea\!\Bigr)
\,/\,
\Bigl(\!\ba{c} \delta+i+j+1/2 \\[1mm] j \ea\!\Bigr)~. \nonumber \\[2mm]
\mbox{}
\end{eqnarray}
Observe that the summation range in (\ref{e125}) is all integer $i\geq0$ such that $j-\delta\leq i\leq j$.

In the corresponding problem of expanding
\beq \label{e128}
{}^2R_{2j}^{0,0}(\rho)(1-\rho^2)^{\delta}=\sum_i\,{}^2F_{ji}^{\delta}\,{}^2R_{2i}^{0,\delta}(\rho)~,
\eq
the required coefficients $F$ are given by
\beq \label{e129}
{}^2F_{ji}^{\delta}=({-}1)^{j-i}\,\frac{\delta+2i+1}{\delta+i+j+1}\,\Bigl(\!\ba{c} \delta \\[1mm] j-i \ea\!\Bigr)\, \Bigl(\!\ba{c} i+j \\[1mm] j \ea\!\Bigr)
\,/\,
\Bigl(\!\ba{c} \delta+i+j \\[1mm] j \ea\!\Bigr)~.
\eq

\section{Expansion of some special radial profiles} \label{sec7}
\mbox{} \\[-9mm]

We present expansion into radial Zernike functions of some special radial profiles as required in Secs.~\ref{sec8} and \ref{sec9} for the analytic construction of all-scale and multi-scale wavelet transforms.

\subsection{Expansion of $\rho^{\beta}(1-\rho^2)^{\alpha}$ as an $R_n^{l,\alpha}$-series} \label{subsec7.1}
\mbox{} \\[-9mm]

We have, see (\ref{e20}),
\beq \label{e130}
B(\rho)=\rho^{\beta}(1-\rho^2)^{\alpha}=\sum_{n=l,l+2,...}\,\frac{b_n^{l,\alpha}}{N_{nl}^{\alpha}}\,R_n^{l,\alpha}(\rho)~,
\eq
where $N_{nl}^{\alpha}$ is given in (\ref{e10}) and
\begin{eqnarray} \label{e131}
b_n^{l,\alpha} & = & \il_0^1\,\rho^{\beta}(1-\rho^2)^{\alpha}\,R_n^{l,\alpha}(\rho)\,\frac{\rho^2\,d\rho}{(1-\rho^2)^{\alpha}} \nonumber \\[3.5mm]
& = & \il_0^1\,\rho^{l+\beta+1}(1-\rho^2)^{\alpha}\,P_p^{(\alpha,l+1/2)}(2\rho^2-1)\,\rho\,d\rho
\end{eqnarray}
with $n=l+2p$ and $l,p=0,1,...\,$. This is the same integral as in (\ref{e41}), except that the power $2l+2k+1$ of $\rho$ is replaced by $l+\beta+1$. Thus replacing $k$ by $\frac12\,(\beta-l)$ in (\ref{e44}), we get
\begin{eqnarray} \label{e134}
b_n^{l,\alpha} & = & \frac12~\frac{\Gamma(\tfrac12\,(\beta-l)+1)}{p!\,\Gamma(\tfrac12\,(\beta-l)-p+1)}~ \frac{\Gamma(p+\alpha+1)\,\Gamma(\tfrac12\,(\beta+l)+3/2)}{\Gamma(\tfrac12\,(\beta+l)+p+\alpha+5/2)} \nonumber \\[3.5mm]
& = & \frac{1}{2\alpha+\beta+l+2p+3}\,\Bigl(\!\ba{c} \dfrac{\beta-l}{2} \\[3mm] p \ea\!\Bigr)
\,/\,
\Bigl(\!\ba{c} \alpha+p+\dfrac{\beta+l+1}{2} \\[3mm] \alpha+p \ea\!\Bigr)~,
\end{eqnarray}
where we use the binomial notation
\beq \label{e135}
\Bigl(\!\ba{c} a \\[1mm] p \ea\!\Bigr)=\frac{a(a-1)\cdot...\cdot(a-p+1)}{p!}~,~~~~~~{\rm any}~a\in\dR\,,~~p=0,1,...~.
\eq

For the 2D case, there is the expansion
\beq \label{e136}
\rho^{\beta}(1-\rho^2)^{\alpha}=\sum_{n=l,l+2,...}\,\frac{{}^2b_n^{l,\alpha}}{{}^2N_{nl}^{\alpha}}\,{}^2R_n^{l,\alpha}(\rho)~,
\eq
with ${}^2N_{nl}^{\alpha}$ given by (\ref{e13}) and
\begin{eqnarray} \label{e137}
{}^2b_n^{l,\alpha} & = & \il_0^1\,\rho^{l+\beta}(1-\rho^2)^{\alpha}\,P_p^{(\alpha,l)}(2\rho^2-1)\,\rho\,d\rho \nonumber \\[3.5mm]
& = & \frac12~\frac{\Gamma(\tfrac12\,(\beta-l)+1)}{p!\,\Gamma(\tfrac12\,(\beta-l)-p+1)}~\frac{\Gamma(p+\alpha+1)\, \Gamma(\tfrac12\,(\beta+l)+1)} {\Gamma(\tfrac12\,(\beta+l)+p+\alpha+2)} \nonumber \\[3.5mm]
& = & \frac{1}{2\alpha+\beta+l+2p+2}\,\Bigl(\!\ba{c} \dfrac{\beta-l}{2} \\[3mm] p \ea\!\Bigr)
\,/\,
\Bigl(\!\ba{c} \alpha+p+\dfrac{\beta+l}{2} \\[3mm] \alpha+p \ea\!\Bigr)~.
\end{eqnarray}

\subsection{Expansion of $\rho^{\beta}(1-\rho^2)^{n+\delta}$ as an $R_n^{l,\delta}$-series in the case that $\beta=l=0,1,...$} \label{subsec7.2}
\mbox{} \\[-9mm]

We have
\beq \label{e138}
\rho^l(1-\rho^2)^{\eta+\delta}=\sum_{n=l,l+2,...}\,\frac{b_n^{l,\delta}}{N_{nl}^{\delta}}\,R_n^{l,\delta}(\rho)~,
\eq
where $N_{nl}^{\delta}$ is given by (\ref{e10}) with $\alpha=\delta$, and
\begin{eqnarray} \label{e139}
b_n^{l,\delta} & = & \il_0^1\,\rho^l(1-\rho^2)^{\eta+\delta}\,R_n^{l,\delta}(\rho)\,\frac{\rho^2\,d\rho}{(1-\rho^2)^{\delta}} \nonumber \\[3.5mm]
& = & \il_0^1\,\rho^{2l+1}(1-\rho^2)^{\eta+\delta}\,P_p^{(\delta,l+1/2)}(2\rho^2-1)\,\rho\,d\rho~.
\end{eqnarray}
Proceeding as in (\ref{e39}--\ref{e43}) with Rodriguez' formula, partial integrations and the Beta-integral, we get now
\begin{eqnarray} \label{e140}
b_n^{l,\delta} & = & \frac12~\frac{({-}1)^p\,\Gamma(\eta+1)}{p!\,\Gamma(\eta-p+1)}~\frac{\Gamma(\delta+\eta+1)\, \Gamma(l+p+3/2)} {\Gamma(\delta+\eta+l+p+5/2)} \nonumber \\[3.5mm]
& = & \frac{({-}1)^p}{2\delta+2\eta+2l+2p+3}\,\Bigl(\!\ba{c} \eta \\[1mm] p \ea\!\Bigr)
\,/\,
\Bigl(\!\ba{c} \delta+\eta+l+p+1/2 \\[1mm] \delta+\eta \ea\!\Bigr)~.
\end{eqnarray}
For the 2D case, for expanding
\beq \label{e141}
\rho^l(1-\rho^2)^{\eta+\delta}=\sum_{n=l,l+2,...}\,\frac{{}^2b_n^{l,\delta}}{{}^2N_{nl}^{\delta}}\,{}^2R_n^{l,\delta}(\rho)~,
\eq
we require ${}^2N_{nl}^{\delta}$ of (\ref{e13}) with $\alpha=\delta$, and
\beq \label{e142}
{}^2b_n^{l,\delta}=\frac{({-}1)^p}{2\delta+2\eta+2l+2p+2}\,\Bigl(\!\ba{c} \eta \\[1mm] p \ea\!\Bigr)
\,/\,
\Bigl(\!\ba{c} \delta+\eta+l+p \\[1mm] \delta+\eta \ea\!\Bigr)~.
\eq

\section{Radial profile for all-scale wavelet transform} \label{sec8}
\mbox{} \\[-9mm]

We consider radial profiles of the form $\rho^{\beta}(1-\rho^2)^{\alpha}$ as a candidate for the radial part $B(\rho)$ of a separable function $F(\vart,\varp,\rho)=A(\vart,\varp)\,B(\rho)$ on the unit ball whose 3D Fourier transform should act as an anisotropic wavelet $\psi$ of an all-scale wavelet. We take $\beta$ an even integer $\geq\,0$ and $\alpha$ an integer $\geq\,0$. The expansion in (\ref{e130}) with coefficients given by (\ref{e137}) is then finite for $l=0,2,...\,$. One can now opt to use these coefficients directly, or, alternatively, to compute them recursively in $l=0,2,...$ as in Sec.~\ref{sec5}, to obtain the 3D Fourier transform of $A_{{\rm Funk}}\,B$ according to Subsec.~\ref{subsec2.5}.

The profile $B$ has its maximum over $\rho\in[0,1]$ at
\beq \label{e143}
\rho_{\max}=\Bigl(\frac{\tfrac12\beta}{\alpha+\tfrac12\beta}\Bigr)^{1/2}~,
\eq
and
\beq \label{e144}
B_{\max}=B(\rho_{\max})=\frac{(\tfrac12\beta)^{\frac12\beta}\,\alpha^{\alpha}}{(\alpha+\tfrac12\beta)^{\alpha+\frac12\beta}}~.
\eq
For somewhat larger values of $\alpha$ and $\beta$, the profile $B$ is rather spiky, and this is considered unfavourable when an all-scale wavelet is desired. We present a simple means to improve $B$ in this respect. We consider, to this end, $B$ as a function of $\rho^2$, and we multiply this function by the second order Taylor expansion of $1/B$ around $\rho^2=\rho_{\max}^2$. Thus, we let
\beq \label{e145}
g(x)=x^{\frac12\beta}(1-x)^{\alpha}~,
\eq
so that $B(\rho)=g(\rho^2)$, and we set
\begin{eqnarray} \label{e146}
h(x) & = & g(x)\,\Bigl(\frac{1}{g_{\max}}+\frac12\,\Bigl(\frac1g\Bigr)''\,(x_{\max})(x-x_{\max})^2\Bigr) \nonumber \\[3.5mm]
& = & \frac{g(x)}{g_{\max}}\,\Bigl(1+\frac12\,\Bigl(\frac1g\Bigr)''\,(x_{\max})\,g_{\max}(x-x_{\max})^2\Bigr)~.
\end{eqnarray}
With $x_{\max}=\rho^2_{\max}$ and $g_{\max}=g(x_{\max})=B_{\max}$, we compute
\beq \label{e147}
\tfrac12\gamma:=\frac12\,\Bigl(\frac1g\Bigr)''\,(x_{\max})\,g_{\max}=\frac{(\alpha+\tfrac12\beta)^3}{\alpha\beta}~.
\eq
Accordingly, we replace $b(\rho)=B(\rho)/B(\rho_{\max})$ by
\beq \label{e148}
c(\rho)=b(\rho)\,\Bigl(1+\frac{(\alpha+\tfrac12\beta)^3}{\alpha\beta}\,(\rho^2-\rho_{\max}^2)^2\Bigr)~.
\eq
We observe that the computation of the expansion of $c(\rho)$ as in (\ref{e130}) is still feasible, since $c(\rho)$ is a linear combination of 3 profiles of the same type and with the same $\alpha$ as $B$ itself.

The above procedure can be repeated, if one wishes, with $h(x)$ in (\ref{e146}) instead of $g(x)$ in (\ref{e145}). It turns out that $h(x)$ is maximal at
\beq \label{e149}
x=\frac{\tfrac14\beta^2+\tfrac12\alpha\beta+2\alpha}{(\alpha+\tfrac12\beta+2)(\alpha+\tfrac12\beta)}~,
\eq
a result that is obtained from a lengthy calculation, rendering this second step analytically feasible.

In Fig.~2a, we show $b(\rho)=B(\rho)/B_{\max}$ where $B(\rho)=\rho^{\beta}(1-\rho^2)^{\alpha}$, $0\leq\rho\leq1$, for $\alpha=6$ and $\beta=2$. This $B$ is maximal at $\rho_{\max}=7^{-1/2}$, with maximal value $B_{\max}=6^6/7^7$. In Fig.~2b, we show $c(\rho)=b(\rho)(1+\frac{1}{12} 7^3(\rho^2-\frac17)^2)$ that flattens $B$ in accordance with (\ref{e148}). This $c$ is maximal at $(19/63)^{1/2}$.
\vspace*{-2cm}
\begin{figure} [!h]
\centerline{\includegraphics[width=0.7\textwidth]{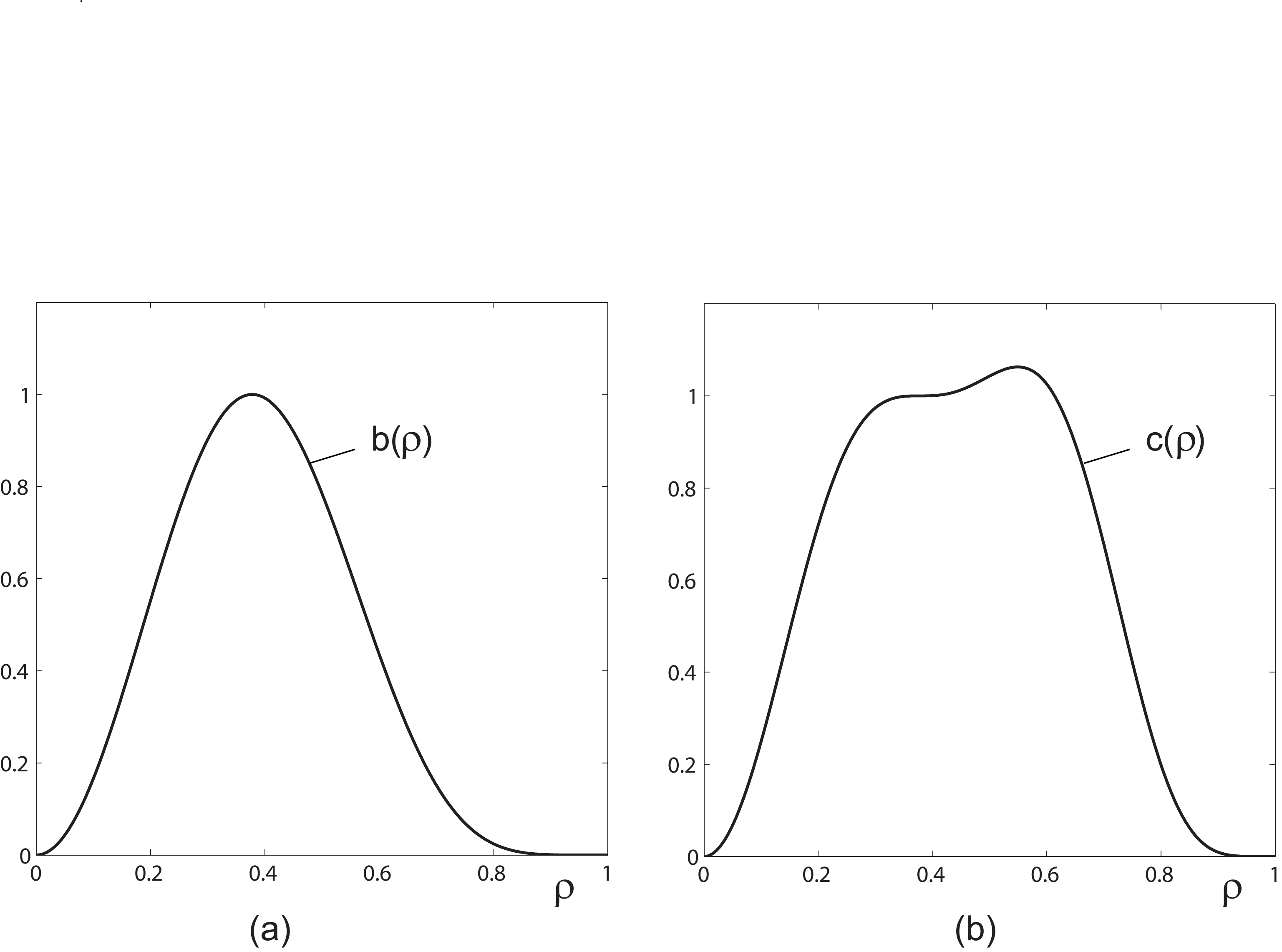}}
\caption{Plot of (a) $b(\rho)=7^76^{-6}\rho^2(1-\rho^2)^6$, $0\leq\rho\leq1$, and of its flattened version (b) $c(\rho)=b(\rho)(1+\frac{1}{12} 7^3(\rho^2-\frac17)^2)$, $0\leq\rho\leq1$, obtained using (\ref{e148}).}
\end{figure}

\newpage
\noindent
\section{Radial profile for multiscale wavelet transform} \label{sec9}
\mbox{} \\[-9mm]

We consider now the setting of Sec.~\ref{sec6} in which we choose
\beq \label{e150}
C(\rho)=(1-\rho^2)^{\eta}~,~~~~~~0\leq\rho\leq1~.
\eq
Thus we consider the scaled-and-truncated standard profiles
\beq \label{e151}
S=S(\rho\,;\,\beta,\delta\,;\,\eta\,;\,\eps)=\rho^{\beta}(1-\eps^2\rho^2)^{\eta}\,(1-\rho^2)^{\delta}~,~~~~~~0\leq\rho\leq1~,
\eq
with fixed truncation parameters $\beta,\delta\geq0$, shape parameter $\eta\geq0$, and scaling parameter $\eps\in[0,1]$. We take $\beta$, $\eta$ and $\delta$ integer with $\beta$ even.

The $C$ in (\ref{e150}) has an explicit expansion in standard 3D radial polynomials $R_{2k}^{0,0}$, $k=0,1,...\,$, which is obtained from Subsec.~\ref{subsec7.2} by choosing $l=\delta=0$ and $\eta=\alpha$ in (\ref{e138}). Thus, for any $\eps$, we can apply the procedure given in Sec.~\ref{sec6} to find the expansion of $S$ in (\ref{e151}) in radial functions $R_{2k}^{0,\delta}$, $k=0,1,...\,$. From this expansion, we can find all expansions of $S$ in radial functions $R_{l+2s}^{l,\delta}$, $s=0,1,...\,$, with $l=0,2,...$ using the recursion in Sec.~\ref{sec5}, see (\ref{e89}). The resulting expansion coefficients can be tested on correctness for the case that $\eps=1$ and $\beta=l$ (even by assumption) by using the explicit expansion of $\rho^l(1-\rho^2)^{\eta+\delta}$ in radial functions $R_n^{l,\alpha}$, $n=l,l+2,...\,$, as given in Subsec.~\ref{subsec7.2}. The whole aggregate of expansion coefficients is then used, as in Subsec.~\ref{subsec2.5}, to find the 3D Fourier transform of $A_{{\rm Funk}}\,S$.

The $S$-profiles (\ref{e151}) vary between the extreme cases
\beq \label{e152}
S_{\eps=1}=\rho^{\beta}(1-\rho^2)^{\eta+\delta}~~{\rm and}~~S_{\eps=0}=\rho^{\beta}(1-\rho^2)^{\delta}~,~~~~~~ 0\leq\rho\leq1~,
\eq
with respective maxima
\beq \label{e153}
\frac{(\tfrac12\beta)^{\frac12\beta}\,(\eta+\delta)^{\eta+\delta}}{(\frac12\beta+\eta+\delta)^{\frac12\beta+\eta+\delta}}~~{\rm and}~~ \frac{(\tfrac12\beta)^{\frac12\beta}\,\delta^{\delta}}{(\tfrac12\beta+\delta)^{\frac12\beta+\delta}}~,
\eq
assumed at
\beq \label{e154}
\rho(1)=\Bigl(\frac{\tfrac12\beta}{\tfrac12\beta+\eta+\delta}\Bigr)^{1/2}~~{\rm and}~~\rho(0)=\Bigl(\frac{\tfrac12\beta}{\tfrac12\beta+\delta}\Bigr)^{1/2}~.
\eq

\subsection{Scaling parameter as a function of ${\rm arg}\,\max\,S(\rho)$} \label{subsec9.1}
\mbox{} \\[-9mm]

The following facts can be established by elementary means. For any $\eps\in[0,1]$, the profile $S$ has a unique maximum at a point $\rho=\rho(\eps)\in[0,1]$. The function $\rho(\eps)$ increases from $\rho(1)$ at $\eps=1$ to $\rho(0)$ at $\eps=0$ when $\eps$ decreases from 1 to 0. Given any value $\overline{\rho}\in[\rho(1),\rho(0)]$, we have $\rho(\eps)=\overline{\rho}$ for
\beq \label{e155}
\eps=\eps(\rho)=\overline{\eps}=\frac{1}{\overline{\rho}}\,\Bigl(\frac{\tfrac12\beta-(\tfrac12\beta+\delta)\,\overline{\rho}^2} {\tfrac12\beta+\eta-(\tfrac12\beta+\eta+\delta)\,\overline{\rho}^2}\Bigr)^{1/2}~.
\eq
The right-hand side of (\ref{e155}) is well-defined for all
\beq \label{e156}
\overline{\rho}\leq\rho(0)\leq\Bigl(\frac{\tfrac12\beta+\eta}{\tfrac12\beta+\eta+\delta}\Bigr)^{1/2}~.
\eq
Allowing $\overline{\rho}<\rho(1)$ in (\ref{e155}) would lead to $\overline{\eps}>1$, which we have excluded. Of course, we can invert (\ref{e155}) so as to obtain $\overline{\rho}$ as a function $\overline{\eps}$, but the resulting formula lacks transparency. For values of $\overline{\rho}$ sufficiently below $\rho(0)$, we have the approximation
\beq \label{e157}
\overline{\eps}\approx\frac{1}{\overline{\rho}}\,\Bigl(\frac{\tfrac12\beta}{\tfrac12\beta+\eta+\delta}\Bigr)^{1/2}=\frac{\rho(1)}{\overline{\rho}}~,
\eq
allowing an interpretation of $\overline{\eps}$ as a true scaling parameter.

\subsection{Design example} \label{subsec9.2}
\mbox{} \\[-9mm]

In designing the parameters $\beta$, $\delta$ and $\eta$ in (\ref{e151}), so that $S$ ranges through a desired set of unimodal profiles when $\eps$ decreases from 1 to 0, the following issues should be taken into account:
\bi{1.0}
\ITEM{1.} the position of the maximum of $S$ ranges between $\rho(1)$ and $\rho(0)$, in accordance with Subsec.~\ref{subsec9.1},
\ITEM{2.} in order that the profiles $S$ accommodate a large range $[\rho(1),\rho(0)]$, while $\overline{\eps}$ of (\ref{e155}) has a credible interpretation as a scaling parameter on a substantial part of that range, one should take $\frac12\beta$ large compared to $\delta$ and $\eta$ large compared to $\frac12\beta+\delta$,
\ITEM{3.} large values of $\delta$ and $\eta$ result into spiky profiles $S$, necessitating a large number of values of the scaling parameter $\eps$ to obtain a sufficient uniform coverage of the range $[\rho(0),\rho(1)]$,
\ITEM{4.} large values of $\delta$ and $\eta$ imply that the various expansion and connection coefficients contain binomials and Pochhammer symbols of high order, and this can be avoided by recursive computation such as was demonstrated for the connection coefficients $C$ at the end of Sec.~\ref{sec5}.
\ei

\begin{figure} [!h]
\centerline{\includegraphics[width=0.6\textwidth]{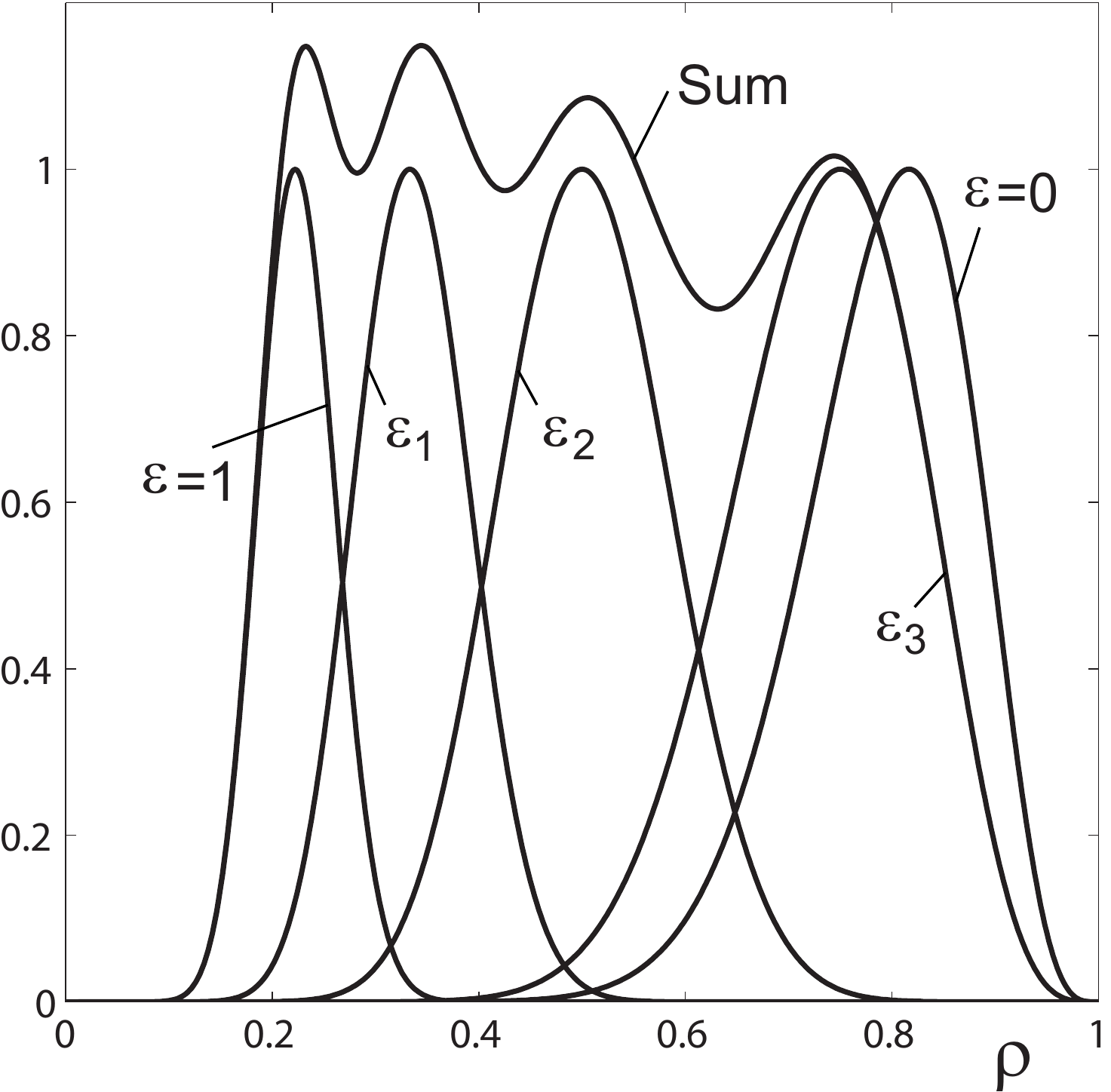}}
\caption{Normalized profiles $S_i(\rho)=\rho^{16}(1-\eps_i^2\rho^2)^{150}\,(1-\rho^2)^4$ for the values $\eps_i=1,(3/7)^{1/2},(8/47)^{1/2},(32/963)^{1/2},0$ of the scaling parameter, yielding maximum positions $\rho_i=2/9,1/3,1/2,3/4,(2/3)^{1/2}$ for $i=0,1,2,3,3.2095=\frac12+{\rm ln}\,3/({\rm ln}(3/2))$, so that $\rho_i=\frac29\,(\frac32)^i$. The sum function is given as $\sum_{i=0}^3\,S_i(\rho)$ and avoids the profile with $\eps_i=0$.}
\end{figure}

\noindent
Some of these issues are illustrated in Fig.~3 that considers the case that $\beta=16$, $\eta=150$ and $\delta=4$ in (\ref{e151}). The value $\delta=4$ provides normally sufficient smoothness and decay of the profiles at the endpoint $\rho=1$. The large value $\beta=16$ provides, in general, more than sufficient decay at the other endpoint $\rho=0$; it has been chose so large to achieve that the $S$-profiles accommodate a $\rho$-range extending all the way to $\rho(0)=(2/3)^{1/2}=0.8165$. Next, the large value $\eta=150$ has been chosen so that the range accommodated by the $S$-profiles starts at a value as low as $\rho(1)=2/9=0.2222$. In Fig.~3, the extreme $S$-profiles with $\eps=0$ and $\eps=1$ occur at the right-most and left-most position, respectively. The three other values of the scaling parameter $\eps$ have been chosen such that the respective maxima occur at $\rho_i=1/3,1/2,3/4$, i.e., at $\frac29\,(\frac32)^i$, $i=1,2,3$. It appears that these three profiles, together with the one with $\eps=1$, provide an adequate coverage of the total range $[\rho=0.18,\rho=0.85]$, as demonstrated by the sum function of these four profiles that is pretty close to being constant on the mentioned range. The values of $\eps$ that achieve the desired positions of the maxima are found using (\ref{e155}). They are given by $\eps_0=1=(2/3)^0$, and by
\beq \label{e158}
\eps_1=\sqrt{3/7}=0.6546\leq(2/3)^1=0.6666~,
\eq
\beq \label{e159}
\eps_2=\sqrt{8/47}=0.4125\leq(2/3)^2=0.4444~,
\eq
\beq \label{e160}
\eps_3=\sqrt{32/963}=0.1882\leq(2/3)^3=0.293
\eq
for the values of the scaling parameters corresponding to $\rho_i=2/9$ and 1/3, 1/2, 3/4, respectively. These values of $\eps_i$ are somewhat, but not dramatically, below $\rho(1)/\rho_i$, see (\ref{e157}). \\ \\
{\bf Acknowledgement.}~~The author wishes to thank M.\ Janssen for providing numerous checks, verifications and implementations of the sometimes complicated formulas in this report. The research leading to this report was carried out under the auspices of R.\ Duits' ERC grant, European Community's 7$^{{\rm th}}$ Framework Programme (FP/2007--2014)/ERC grant agreement No.\ 335555.

\end{document}